\documentclass{article}
\usepackage{graphicx}
\usepackage{graphicx}
\usepackage{array,graphicx}
\usepackage{epstopdf}
\usepackage{algorithm}
\usepackage{algpseudocode}
\usepackage{amsopn}
\usepackage{amsfonts}
\usepackage{amsmath}
\usepackage{amssymb}
\usepackage{mathtools}
\usepackage{color}
\usepackage{siunitx}
\usepackage{url,latexsym,enumerate,graphics,times,commath}
\usepackage{tikz}
\usepackage{pgfplots}
\usepackage{subcaption}
\usepackage[hidelinks]{hyperref}
\newlength\mylen
\newcommand*\Lap{\mathrm{\Delta}}

\ifpdf
\DeclareGraphicsExtensions{.pdf,.png,.jpg}
\else
\DeclareGraphicsExtensions{.eps}
\fi

\newcommand*\samethanks[1][\value{footnote}]{\footnotemark[#1]}
\newcommand{\specificthanks}[1]{\@fnsymbol{#1}}
\ifpdf
\hypersetup{
	pdftitle={A meshfree collocation scheme for surface differential operators on point clouds},
	pdfauthor={Abhinav Singh, Ale Foggia, Pietro Incardona, and  Ivo F.~Sbalzarini}
}
\fi
\begin{document}
	\title{A meshfree collocation scheme for surface differential operators on point clouds}
	\date{}
	\author{Abhinav Singh\thanks{Technische Universit\"{a}t Dresden, Faculty of Computer Science, Dresden, Germany. \newline \& Max Planck Institute of Molecular Cell Biology and Genetics, Dresden, Germany.\newline \& Center for Systems Biology Dresden, Dresden, Germany.}
		\and Alejandra Foggia\samethanks[1]
		\and Pietro Incardona\samethanks[1]
		\and Ivo F.~Sbalzarini\samethanks[1] \thanks{Cluster of Excellence Physics of Life, TU Dresden, Dresden, Germany.}
	}			\maketitle
	{\let\thefootnote\relax\footnotetext{This work was funded by the German Research Foundation (DFG, Deutsche Forschungsgemeinschaft) under grants FOR-3013 (``Vector- and tensor-valued surface PDEs''), GRK-1907 
			(``RoSI: role-based software infrastructures''), and SB-350008342 (``OpenPME''), and by the Federal Ministry of Education and Research (Bundesministerium f\"{u}r Bildung und Forschung, BMBF) under funding code 031L0160 (project ``SPlaT-DM -- computer simulation platform for topology-driven morphogenesis'').}}
			\begin{abstract}
We present a meshfree collocation scheme to discretize intrinsic surface differential operators over scalar fields on smooth curved surfaces with given normal vectors and a non-intersecting tubular neighborhood. The method is based on Discretization-Corrected Particle Strength Exchange (DC-PSE), which generalizes finite difference methods to meshfree point clouds. The proposed Surface DC-PSE method is derived from an embedding theorem, but we analytically reduce the operator kernels along surface normals to obtain a purely intrinsic computational scheme over surface point clouds. We benchmark Surface DC-PSE by discretizing the Laplace-Beltrami operator on a circle and a sphere, and we present convergence results for both explicit and implicit solvers. We then showcase the algorithm on the problem of computing Gauss and mean curvature of an ellipsoid and of the Stanford Bunny by approximating the intrinsic divergence of the normal vector field. Finally, we compare Surface DC-PSE with Surface Finite Elements (SFEM) and Diffuse-Interface Finite Elements (DI FEM) in a validation case.
			\end{abstract}
			
			\section{Introduction}
Partial Differential Equations (PDEs) on curved surfaces and differentiable manifolds are an important tool in understanding and studying physical phenomena such as surface flows~\cite{voigt_fluid_2019,nitschke_liquid_2019} and active morphogenesis~\cite{mietke_self-organized_2019}. Analytically solving intrinsic PDEs in curved surfaces, however, quickly becomes impossible for nonlinear PDEs or for surfaces that do not possess a global parameterization. Therefore, numerical methods for solving intrinsic PDEs on curved surfaces are important, and a wide variety of both embedded and embedding-free schemes have been developed to consistently discretize intrinsic differential operators over scalar fields on surfaces.

Embedding-free methods require a (at least local) parameterization of the surface in order to discretize the differential operators via coordinate charts or a local basis of the manifold~\cite{wang_discretizing_2013}. This includes methods based on local moving frames \cite{chun2012method}, a concept originally developed in continuous group theory, where the surface geometry is locally represented by intrinsic orthonormal bases. This has been used to solve surface PDEs over meshfree point clouds using Moving Least Squares (MLS) approximations~\cite{liang_solving_2013}. The concept of local moving frames  has also been combined with discontinuous Galerkin discretization, e.g., to solve shallow-water equations on arbitrary rotating surfaces~\cite{chun2017method}. Other embedding-free Finite Element Methods (FEM) include Intrinsic Surface FEM (ISFEM), which discretizes differential operators on a triangulation of the surface \cite{bachini_intrinsic_2021,grande_space-time_2014}, and methods based on Discrete Exterior Calculus (DEC) \cite{nitschke_discrete_2017}.

Embedding methods discretize the surface problem in an embedding space of co-dimension 1 and use projections to restrict the differential operators computed in the embedding space to the surface manifold. This includes  methods that use explicit tracer points to represent the surface, but interpolate to an embedding mesh to evaluate differential operators~\cite{Leung:2009b}, diffuse-interface methods based on phase-field representations of the surface~\cite{nestler_finite_2019}, embedding FEM such as TraceFEM~\cite{olshanskii_trace_2017} and diffuse-interface FEM~\cite{Lehrenfeld2017,diffuse_interface_2006}, narrow-band level-set methods based on orthogonal extension of the surface quantities~\cite{cottet_semi-implicit_2016,bergdorf_lagrangian_2010}, level-set methods based on the closest-point transform~\cite{ruuth_simple_2008,macdonald_solving_2011}, and volume-of-fluid methods for surface PDE problems~\cite{James:2004}.

While each of these methods has its specific strengths, embedding methods usually generalize better to complex-shaped or arbitrary surfaces~\cite{ruuth_simple_2008}. However, they tend to have higher computational cost, because computations are done in the higher-dimensional embedding space and additional extension (for level sets), right-hand-side evaluation (for phase fields), or interpolation (for closest-point transforms) steps are required, albeit specific optimizations are available, e.g., for level sets~\cite{Peng:1999}. Embedding-free methods are usually more accurate, because they avoid the interpolation and projection errors arising when the discretization of the embedding space does not trace the surface exactly, but they tend to be more difficult to implement and harder to generalize to complex-shaped or moving/deforming surfaces.

Here, we present a meshfree collocation method for PDEs on smooth and orientable curved surfaces with non-intersecting tubular neighborhood. The method combines elements from embedding and embedding-free approaches. It is {\em algorithmically} embedding-free in the sense that surface quantities are represented on tracer points that are contained in the surface. This also discretizes and represents the surface itself as a point cloud. But the method is {\em mathematically} related to embedding approaches, since the stencils used to approximate differential operators at the surface points are computed in the embedding space by a reduction operation along the local normal vector, which needs to be known or computed. Intuitively, this projects the discrete operators, rather than projecting the flux vectors as typically done in embedding methods. The resulting method therefore shares properties of moving frame approaches, such as the low dimensionality (and hence low computational cost) and the meshfree character~\cite{chun2012method,liang_solving_2013}. It combines these with properties of embedding methods, such as their flexibility in generalizing to complex surfaces~\cite{ruuth_simple_2008}, and their ability to compute extrinsic differential-geometric quantities.

Our method is based on the Discretization-Corrected Particle Strength Exchange (DC-PSE) collocation scheme for arbitrary (surface) point clouds.
DC-PSE is related to Generalized Finite Difference Methods (GFDM) \cite{suchde_meshfree_2018} and to MLS \cite{schrader_discretization_2010}.
Given the local surface normal $\boldsymbol{n}$, we derive intrinsic discrete operators by first creating an embedding narrow-band and placing collocation points along the normal from each surface point. We then determine the regular DC-PSE operator kernels in the embedding space. These kernels are subsequently reduced under the condition of orthogonal extension $\nabla f \cdot \boldsymbol{n}=0$ for any (sufficiently) differentiable scalar field $f(\boldsymbol{x})\in\mathbb{R}$ to derive intrinsic kernels at the surface points $\boldsymbol{x}$. This is possible due to the kernel nature of DC-PSE, and it preserves the information from the embedding space in a scheme that only requires computation over surface points.

This paper is organized as follows: Section \ref{sec:DCPSE} recollects the DC-PSE method for convenience and introduces the notation. In Section \ref{sec:SDCPSE}, we describe the Surface DC-PSE scheme for numerically consistent discretization of surface differential operators. We present validation and convergence result in Section \ref{sec:Result} and conclude in Section \ref{sec:Conclusion}.

\section{Discretization-Corrected Particle Strength Exchange (DC-PSE)}\label{sec:DCPSE}
DC-PSE is a numerical method for discretizing differential operators on irregular distributions of collocation points~\cite{schrader_discretization_2010}. The method was originally derived as an improvement over the classic Particle Strength Exchange (PSE)~\cite{eldredge_general_2002} scheme, reducing its quadrature error on irregularly distributed collocation points, but mathematically amounts to a generalization of finite differences~\cite{schrader_discretization_2010}. The PSE/DC-PSE class of collocation methods uses mollification with a symmetric smoothing kernel $\eta_\epsilon (\cdot )$ to approximate (sufficiently smooth) continuous functions $f(\boldsymbol{x})\in\mathbb{R}$,\, $\boldsymbol{x}\in\mathrm{\Omega}\subseteq \mathbb{R}^d$,
\begin{equation}
	f(\boldsymbol{x}_{p}) \approx f_{\epsilon}(\boldsymbol{x}_{p})=\int_{\mathrm{\Omega}} f(\boldsymbol{x})\, \eta_{\epsilon}(\boldsymbol{x}_{p}-\boldsymbol{x}) \,\mathrm{d} \boldsymbol{x},
	\label{eq:moll}
\end{equation}
where $f_\epsilon(\boldsymbol{x}_{p})$ is a regularized approximation of the function $f$ at location $\boldsymbol{x}_{p}\in \mathrm{\Omega}$ of collocation point $p$. The scalar $\epsilon$ is the smoothing length (or the kernel width) of the mollification.
Linear differential operators in $\mathbb{R}^d$,
\begin{equation}
	\boldsymbol{D}^{\boldsymbol{\alpha}} = \frac{\partial^{| \boldsymbol{\alpha}|}}{\partial x_1^{\alpha _1} \partial x_2^{\alpha _2} \cdots \partial x_d^{\alpha _d}}\, ,
\end{equation}
defined by the multi-index $\boldsymbol{\alpha}=(\alpha _1, \ldots , \alpha _d)\in\mathbb{Z}^d$ with $|\boldsymbol{\alpha}| = \sum_{i=1}^d \alpha_i$ are approximated by Taylor series expansion to find a discrete operator
\begin{equation}
	\boldsymbol{Q}^{\boldsymbol{\alpha}} f(\boldsymbol{x}_{p})=\boldsymbol{D}^{\boldsymbol{\alpha}} f(\boldsymbol{x}_{p})+{O}\!\left(h(\boldsymbol{x}_{p})^{r}\right)
\end{equation}
at collocation point $\boldsymbol{x}_p$. The order of approximation $r$ depends on the kernel $\eta_\epsilon$ used in Eq.~(\ref{eq:moll}), and $h(\boldsymbol{x}_{p})$ is the average distance between collocation point $p$ and its neighbors within the kernel support.
We use the arithmetic mean of the $L_1$-distances to compute $h$, but since all norms are equivalent, the convergence order (but not the actual error magnitude) is independent of the choice of average.
The Taylor expansion yields integral constraints (also known as {\em continuous moment conditions}), which the kernel $\eta_\epsilon$ needs to fulfill in order to reach a certain convergence order $r$~\cite{eldredge_general_2002}.

\begin{figure}[t]
	\setlength{\tabcolsep}{1pt}
	\centering
	\begin{tabular}{cc}
		\subfloat[]{\includegraphics[width=0.49\textwidth]{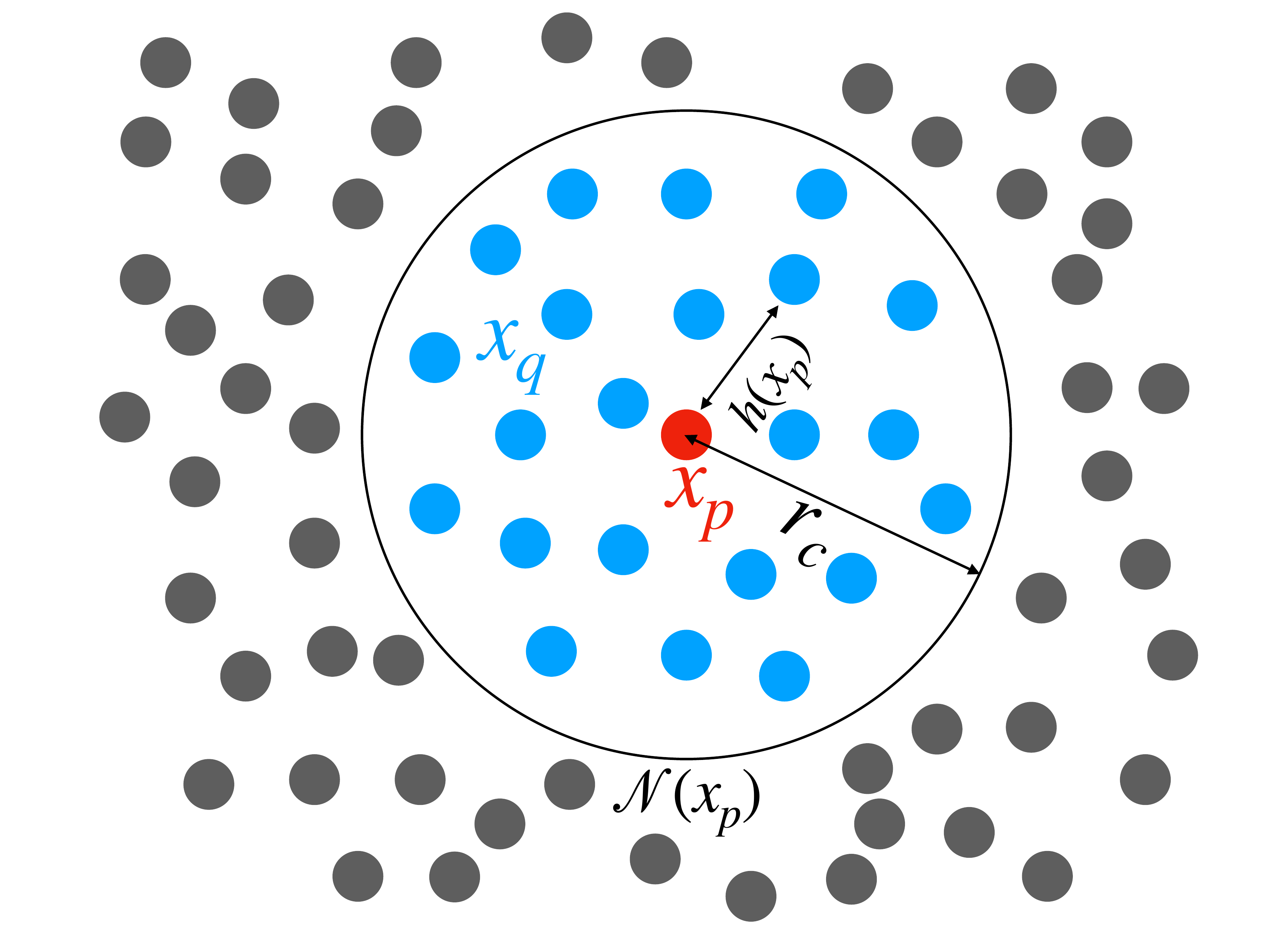}} & \subfloat[]{\includegraphics[width=0.49\textwidth]{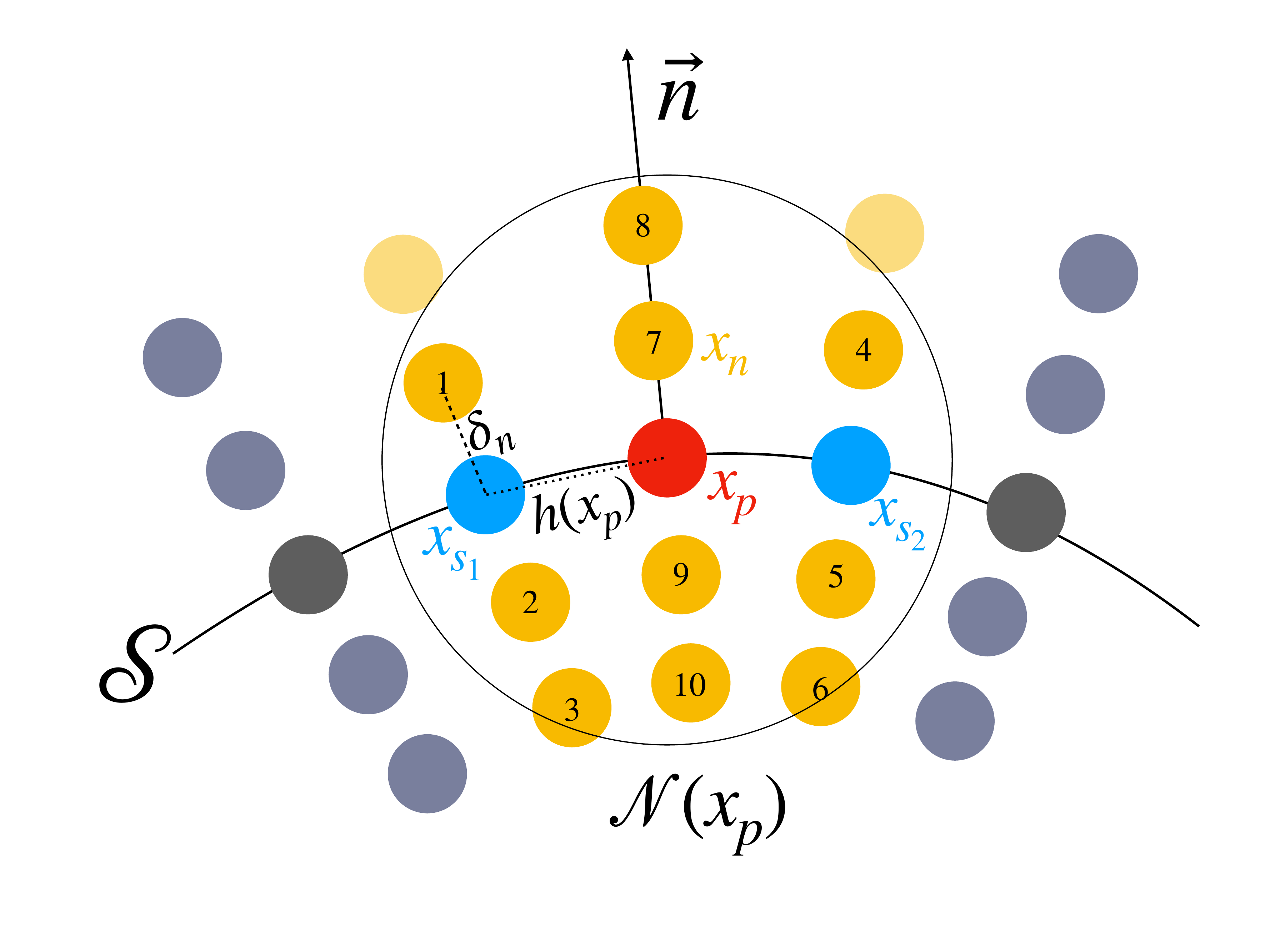}}
	\end{tabular}
   \caption{\textbf{(a) Illustration of the DC-PSE method}. The collocation points $\boldsymbol{x}_{q}$ (blue) within the symmetric operator support $\mathcal{N}(\boldsymbol{x}_p)$ of radius $r_c$ around the center point $\boldsymbol{x}_p$ (red) are used to approximate the differential operator at $\boldsymbol{x}_p$. The average distance between points in the operator support, $h(\boldsymbol{x}_p)$ defines the accuracy of the approximation. \textbf{(b) Illustration of the Surface DC-PSE method.} The intrinsic differential operator at a point $\boldsymbol{x}_p$ (red) on the surface $\mathcal{S}$ is evaluated over neighboring surface points $\boldsymbol{x}_{s}$ (blue). The operator kernel is constructed in the tubular neighborhood of radius $r_c$ (circle) with $2N_n$ points $\boldsymbol{x}_n$ (yellow) replicated along the local surface normal $\boldsymbol{n}$ at distances $\delta n$. Point labels are used in the main text for the illustrative example.}
	\label{fig:sdcpse}
\end{figure}

DC-PSE uses different kernels $\eta_{\epsilon}^{p}(\cdot ,\cdot )$ for different collocation points $p$ and directly acts on a given quadrature of Eq.~(\ref{eq:moll}) with collocation points $\boldsymbol{x}_q \in \mathrm{\Omega}$, resulting in the discrete operator:
\begin{equation}
	\boldsymbol{Q}_{h}^{\boldsymbol{\alpha}} f(\boldsymbol{x}_p)=\frac{1}{\epsilon (\boldsymbol{x}_p)^{|\boldsymbol{\alpha}|}} \sum_{\boldsymbol{x}_q \in \mathcal{N}(\boldsymbol{x}_p)}\!\!\! \left(f(\boldsymbol{x}_{q}) \pm f(\boldsymbol{x}_p)\right) \eta_{\epsilon}^{p}(\boldsymbol{x}_p,\boldsymbol{x}_{q}),
	\label{eq:dcpseEval}
\end{equation}
where $\mathcal{N}(\boldsymbol{x}_{p})$ are all collocation points in the neighborhood (of a certain radius $r_c$ defined by the kernel width) around point $\boldsymbol{x}_{p}$, as illustrated in Fig.~\ref{fig:sdcpse}a. The positive sign in the parenthesis is used for odd $|\boldsymbol{\alpha}|$, the negative sign for even $|\boldsymbol{\alpha}|$. This renders the operator conservative on symmetric collocation point distributions, i.e., when $\eta_{\epsilon}^{p}(\boldsymbol{x}_p,\boldsymbol{x}_{q}) = \eta_{\epsilon}^{q}(\boldsymbol{x}_q,\boldsymbol{x}_{p})$. In DC-PSE, the kernels $\eta_{\epsilon}^{p}$ are thus not determined from continuous moment conditions, as in PSE, but directly from the {\em discrete moment conditions} that result from substituting Eq.~(\ref{eq:dcpseEval}) into the quadrature of Eq.~(\ref{eq:moll})~\cite{schrader_discretization_2010} for a given set $\{\boldsymbol{x}_q\}_{q=1}^N$. This adapts the kernels to the specific distribution of collocation points (hence the name ``discretization-corrected'') and avoids the quadrature error of PSE~\cite{eldredge_general_2002}, leading to a scheme that is consistent with order $r$ on almost\footnote{The collocation point distribution must not be degenerate in the sense that the Vandermonde matrix of the kernel system must have full rank~\cite{bourantas_using_2016}. A trivial example: placing all points along a line and then asking for an approximation of the derivative in the perpendicular direction cannot work.} arbitrary collocation point sets. This means that at each collocation point, a potentially different kernel is used for the same differential operator if the neighboring collocation points within the kernel support are distributed differently.
Evaluating such a kernel at the locations of the collocation points yields a generalized finite-difference stencil, which reduces to the classic compact finite differences on regular grid arrangements of points~\cite{schrader_discretization_2010}.

DC-PSE kernels are determined at runtime by solving a small system of linear equations for each collocation point, resulting from the discrete moment conditions in its kernel neighborhood.
For this, one can choose the function space such that the kernels are compact and symmetric. A frequent choice are polynomials windowed by truncated exponentials~\cite{schrader_choosing_2012}
\begin{equation}
	\eta^{p}_\epsilon(\boldsymbol{x}_p, \boldsymbol{x}_q)=\eta^{p}_\epsilon\left(\frac{\boldsymbol{x}_p- \boldsymbol{x}_q}{\epsilon (\boldsymbol{x}_p)}\right) :=
	\left(\sum_{|\gamma|=\beta_{\min }}^{|\boldsymbol{\alpha}|+r-1} a_{\gamma}(\boldsymbol{x}_p) \left(\frac{\boldsymbol{x}_p-\boldsymbol{x}_q}{\epsilon(\boldsymbol{x}_p)}\right)^{\!\!\gamma}\right) \textrm{e}^{-\left|\frac{\boldsymbol{x}_p-\boldsymbol{x}_q}{\epsilon(\boldsymbol{x}_p)}\right|^{2}}
	\label{eq:dcpseker}
\end{equation}
of finite radius $r_c$. The polynomial coefficients $a_{\gamma}$ are determined for a given $\boldsymbol{\alpha}$ and given collocation points $\boldsymbol{x}_q\in\mathcal{N}(\boldsymbol{x}_p)$, such that the following discrete moment conditions are satisfied:
\begin{equation}
	Z^{\boldsymbol{\beta}}_h=\left\{\begin{array}{ll}
		(-1)^{|\boldsymbol{\alpha}|} \boldsymbol{\alpha} !, \quad \boldsymbol{\beta}=\boldsymbol{\alpha}\\
		0, \qquad \boldsymbol{\beta} \neq \boldsymbol{\alpha}, \quad \beta_{\min} \leqslant|\boldsymbol{\beta}| \leqslant|\boldsymbol{\alpha}|+r-1, ~\beta_{\min}=\left\{\begin{array}{ll}0, ~|\boldsymbol{\alpha}| \text{ odd}\\ 1, ~|\boldsymbol{\alpha}| \text{ even} \end{array} \right.
		\\
		<\infty , \qquad |\boldsymbol{\beta}|=|\boldsymbol{\alpha}|+r\end{array}\right.
\end{equation}
where
\begin{equation}
	Z_{h}^{\boldsymbol{\beta}}(\boldsymbol{x})=\frac{1}{\epsilon(\boldsymbol{x}_p)^{d}} \sum_{\boldsymbol{x}_q \in \mathcal{N}(\boldsymbol{x}_p)} \!\!\frac{(\boldsymbol{x}_p-\boldsymbol{x}_{q})^{\boldsymbol{\beta}}}{\epsilon(\boldsymbol{x}_p)^{|\boldsymbol{\beta}|}} \,\eta^{p}_\epsilon \! \left(\frac{\boldsymbol{x}_p-\boldsymbol{x}_{q}}{\epsilon(\boldsymbol{x}_p)}\right)
\end{equation}
is the discrete moment of order $\boldsymbol{\beta}$ of the kernel $\eta_{\epsilon}^{\boldsymbol{\alpha}}$, and $\beta_{\min}$ is the parity of $|\boldsymbol{\alpha}|$, because the zeroth moment $Z_h^\bold{0}$ vanishes for even operators.
Under these conditions, DC-PSE is consistent with order $r$ as long as
\begin{equation}
	\frac{h(\boldsymbol{x}_{p})}{\epsilon(\boldsymbol{x}_{p})} \in {O}(1),
\end{equation}
i.e., the kernel width $\epsilon$ scales proportionally with the average inter-point distance $h$ around $\boldsymbol{x}_p$ \cite{schrader_discretization_2010}.

\section{Surface DC-PSE}\label{sec:SDCPSE}
We generalize DC-PSE to surface differential operators based on the following classic result~\cite{ruuth_simple_2008,Marz:2012}:
Let $\mathcal{S}\subset \mathbb{R}^d$ be a differentiable manifold that possesses a tubular neighborhood $T$ and is orientable\footnote{Every boundary-less smooth surface embedded in $\mathbb{R}^d$ has a tubular neighborhood, and the orientability condition is not restrictive when considered locally \cite{Marz:2012}.}
and $f:\mathcal{S}\rightarrow\mathbb{R}$. Define $F:T\rightarrow\mathbb{R}$, such that the restriction $\left. F\right|_\mathcal{S}=f$, and $F$ is constant along the  normal direction $\boldsymbol{n}$ of $\mathcal{S}$, i.e., $\nabla F \cdot \boldsymbol{n}=0$. Then, on the surface $\mathcal{S}$,
\begin{equation}
	\nabla_{\!\!\mathcal{S}} f=( \nabla F)|_\mathcal{S},
	\label{eq:surfaceIdentity}
\end{equation}
where $ \nabla_{\!\mathcal{S}} f$ is the intrinsic surface gradient.
A similar result is true for the intrinsic divergence operator ($\nabla_{\!\mathcal{S}}\,\cdot$ ) and for tangential vector field that is extended by constant extension to all surfaces displaced along the normal of $\mathcal{S}$~\cite{ruuth_simple_2008,Marz:2012}.

Given this result, it is straightforward to see the advantages of a meshfree discretization: it allows for conforming discretization of the surface and for {\em exact} constant orthogonal extension by simply copying points along the normal. This creates an embedding narrow-band of {\em exact} closest-point function values within the tubular neighborhood $T$ without a need for interpolation. If $T$ is non-intersecting with a radius of at least $r_c$ everywhere, the conditions of the result in Eq.~(\ref{eq:surfaceIdentity}) are satisfied in the constructed embedding. Analogous to the closest-point method \cite{ruuth_simple_2008}, one can then discretize differential operators in the embedding space. Due to the additive kernel nature of DC-PSE, the discrete operators in the embedding space can be reduced to only the surface points.

This reduction becomes clear from the formulation of the DC-PSE method. Indeed, we realize that the constant normal extension can be made internal to the operator evaluation by accumulating the kernel {\em coefficients} along the normals. To see this, consider the DC-PSE operator in Eq.~(\ref{eq:dcpseEval}) in the embedding space. The neighborhood $\mathcal{N}$ for the summation contains both surface points $\boldsymbol{x}_{s}$ and  normally extended points $\boldsymbol{x}_{n}$, as shown in Fig.~\ref{fig:sdcpse}b. Because the $f(\boldsymbol{x}_{n})$ are identical copies of the values of the respective surface points, we note that the pre-factors $\left(f(\boldsymbol{x}_{s}) \pm f(\boldsymbol{x}_{p})\right)$ in the summation of Eq.~(\ref{eq:dcpseEval}) are the same for all extended normal points and the corresponding surface point $\boldsymbol{x}_{s}$. Hence, for each given pair of a center point $\boldsymbol{x}_p$ and another surface point $\boldsymbol{x}_s$, the interactions with the corresponding normally extended points can be factored out from the kernel summation:
\begin{multline}
	\frac{f(\boldsymbol{x}_{s}) \pm f(\boldsymbol{x}_{p})}{\epsilon(\boldsymbol{x}_{p})^{|\alpha|}} \sum_{\boldsymbol{x}_{q}=\{\boldsymbol{x}_s, \boldsymbol{x}_n :  (\boldsymbol{x}_{n}-\boldsymbol{x}_{s}) || \boldsymbol{n}(\boldsymbol{x}_s) \} }
	\eta_\epsilon^{p}\left(\frac{\boldsymbol{x}_{p}-\boldsymbol{x}_{q}}{\epsilon\left(\boldsymbol{x}_{p}\right)}\right) \\
	= \frac{f(\boldsymbol{x}_{s}) \pm f(\boldsymbol{x}_{p})}{\epsilon(\boldsymbol{x}_{p})^{|\alpha|}}\eta_\mathcal{S}(\boldsymbol{x}_p,\boldsymbol{x}_s),
	\label{eq:sdcpseker}
\end{multline}
defining the surface kernels $\eta_\mathcal{S}(\boldsymbol{x}_p,\boldsymbol{x}_s)$. These can be evaluated over only the surface points $\boldsymbol{x}_s = \mathcal{N}_\mathcal{S}(\boldsymbol{x}_{p})$ in the in-surface neighborhood $\mathcal{N}_\mathcal{S}(\boldsymbol{x}_{p})$ around the surface point $\boldsymbol{x}_{p}$, see Fig.~\ref{fig:sdcpse}b, yielding the Surface DC-PSE operator:
\begin{equation}
	\boldsymbol{Q}_\mathcal{S}^{\boldsymbol{\alpha}} f(\boldsymbol{x}_{p})=\frac{1}{\epsilon(\boldsymbol{x}_{p})^{|\alpha|}} \sum_{\boldsymbol{x}_{s} \in \mathcal{N}_\mathcal{S}(\boldsymbol{x}_{p})}\left(f(\boldsymbol{x}_{s}) \pm f(\boldsymbol{x}_{p})\right) \eta_\mathcal{S}(\boldsymbol{x}_p,\boldsymbol{x}_{s}) .
	\label{eq:sdcpseEval}
\end{equation}
Importantly, the surface kernels $\eta_\mathcal{S}(\boldsymbol{x}_p,\boldsymbol{x}_{s})$, summed over all orthogonally extended points, can directly be computed when determining the kernel weights and without explicitly creating or storing the normally extended points $\boldsymbol{x}_n$.

Evaluating a Surface DC-PSE operator involves only the neighboring points on the surface and requires no narrow band or normally extended grid, even though the construction of the operators uses an embedding. This leads to a corresponding reduction in computational complexity for operator evaluation, as computations are only performed on a $(d-1)$-dimensional surface embedded in $d$-dimensional space. In comparison, the cost of operator evaluation for embedding methods such as the closest-point method is $O(k(d-1))$, where $k>1$ is the narrow-band width, which scales proportionally with the order of convergence.

\subsection{Surface DC-PSE kernel construction}

Surface DC-PSE requires two algorithms that are not part of the standard, flat-space DC-PSE method:
an algorithm to create the intrinsic neighborhood of a surface point $p$, and an algorithm to determine the surface kernels $\eta_{\mathcal{S}}$ at a surface point $p$. We follow the example of Ref.~\cite{bourantas_using_2016} and use explicit component notation and a concrete example in order to directly relate to implementations in computer code.

Figure~\ref{fig:sdcpse}b illustrates a piece of the tubular neighborhood of radius $r_c$ (circle) of a curved surface $\mathcal{S}$ embedded in $\mathbb{R}^2$. The red point $p$ at position $\boldsymbol{x}_p$ on the surface is the ``center'' collocation point at which we derive the discrete Surface DC-PSE operator $\boldsymbol{Q}_\mathcal{S}^{\boldsymbol{\alpha}} f(\boldsymbol{x}_{p})$ for a scalar surface field $f$. Points in light blue are surface points within the embedding-space neighborhood of radius $r_c$ (circle) of $\boldsymbol{x}_p$, and the yellow points are the orthogonal extensions $\boldsymbol{x}_n$. By default, the spatial separation $\delta n$ between adjacent orthogonal extensions of the surface point $p$ is the arithmetic mean of the distances between $p$ and the other surface points within the kernel support, measured in the embedding space. This favors isotropic-resolution neighborhoods and, thus, low condition numbers of the DC-PSE kernel system matrix. The number of orthogonally extended points should be $N_n \approx \mathrm{round}(r_c/\delta n)$ to either side of the surface, which is the default. Only surface points are actually allocated and stored.

In order to determine the DC-PSE kernel $\eta_{\epsilon}^{p}$ in the embedding space, the distances between $\boldsymbol{x}_p$ and all collocation points in its embedding-space neighborhood $\mathcal{N}$ are required. In the example of Fig.~\ref{fig:sdcpse}b, the neighborhood (circle) includes the surface points $\mathcal{N}_\mathcal{S} = \{s_1,s_2\}$ and the normally extended points $\{n_i\}_{i=1}^{10}$. The two pale-yellow points are not part of the neighborhood. Algorithm \ref{alg:embedding_neighborhood} constructs the neighbor set along with the corresponding distances.
Surface normals at a given point are indexed by the point index for better readability, i.e., $\boldsymbol{n}(\boldsymbol{x}_p):=\boldsymbol{n}_p$. In the example of the figure, this results in the output $$\mathcal{N}_{\text{dist}}(\boldsymbol{x}_p)=
   [[\boldsymbol{d}_{s_1},\boldsymbol{d}_{n_1},\boldsymbol{d}_{n_2},\boldsymbol{d}_{n_3}], [\boldsymbol{d}_{s_2},\boldsymbol{d}_{n_4},\boldsymbol{d}_{n_5},\boldsymbol{d}_{n_6}],[\boldsymbol{d}_{n_7},\boldsymbol{d}_{n_8},\boldsymbol{d}_{n_9},\boldsymbol{d}_{n_{10}}]],$$
   where $\boldsymbol{d}_q$ is the distance between the collocation points $p$ and $q$ in the embedding space in units of $\epsilon_p:=\epsilon (\boldsymbol{x}_p)$.

\begin{algorithm}[t]
	\textbf{Input:}
	\begin{enumerate}
      \item Point set $\mathsf{P}$ on the surface $\mathcal{S}$
		\item Cutoff radius for the operator support $r_c$
		\item Indices $\mathcal{N}_{\mathcal{S}}$ of surface points in the neighborhood of $p$
		\item Optional: spacing $\delta n$ between the normally extended points. Default: average embedding-space distance between surface points
		\item Optional: Number of normal copies of each surface point to be used during operator construction $N_n$ (symmetric to either sides of the surface). Default: $N_n=\mathrm{round}(r_c/\delta n)$
	\end{enumerate}
	\textbf{Output:}
   List of distances between point $p$ and all surface $s_i$ and normal $n_i$ points in its neighborhood $\mathcal{N}(\boldsymbol{x}_p)$:  $\mathcal{N}_{\text{dist}}(\boldsymbol{x}_p$)
   \\
	\begin{algorithmic}[1]
		\Require $|\boldsymbol{n}_{s_i}| = |\boldsymbol{n}_p| = 1$
      \State  $\mathcal{N}_{\text{dist}} = [\,]$, $k=0$
		\ForAll{$s_i \in \mathcal{N}_{\mathcal{S}}$}
         \State $\mathcal{N}_{\text{dist}}.\texttt{append}([\,])$
			\For{$i \in [-N_n,N_n]$}
				\State $\boldsymbol{d}_{n_i} = \boldsymbol{x}_p -\boldsymbol{x}_{s_i} - i  \delta n \cdot \boldsymbol{n}_{s_i}$
				\If{$|\boldsymbol{d}_{n_i}| \le r_c$}
            \State $\mathcal{N}_{\text{dist}}[k].\texttt{append} (\boldsymbol{d}_{n_i}/\epsilon_p)$ \Comment{Add to set of the corresponding $s_i$.}
				\EndIf
			\EndFor
			\State $k  \mathrel{+}= 1$
		\EndFor
		\Statex
      \State $\mathcal{N}_{\text{dist}}.\texttt{append}([\,])$
		\For{$i \in [-N_n,N_n] \backslash \{0\}$} \Comment{Normal points to $p$.}
			\State $\boldsymbol{d}_{n_i} = - i  \delta n \cdot \boldsymbol{n}_p$
			\If{$|\boldsymbol{d}_{n_i}| \le r_c$}
         \State $\mathcal{N}_{\text{dist}}[k].\texttt{append}(\boldsymbol{d}_{n_i}/\epsilon_p)$ \Comment{Create a new set containing all points along $\boldsymbol{n}_p$.}
			\EndIf
		\EndFor
	\end{algorithmic}
	\caption{Surface DC-PSE: construction of neighborhood of a point $p$.}
	\label{alg:embedding_neighborhood}
\end{algorithm}

Using this neighborhood data structure, the embedding-space DC-PSE operator at point $p$ in the example of the figure reads:
\begin{align}
	\nonumber
	\boldsymbol{Q}_\mathcal{S}^{\boldsymbol{\alpha}} f(\boldsymbol{x}_{p}) = &
	\frac{f(\boldsymbol{x}_{s_1}) \pm f(\boldsymbol{x}_{p})}{\epsilon(\boldsymbol{x}_{p})} \left( \eta_{\epsilon}^{p}(\boldsymbol{d}_{s_1}) +  \eta_{\epsilon}^{p}(\boldsymbol{d}_{n_1}) + \eta_{\epsilon}^{p}(\boldsymbol{d}_{n_2}) + \eta_{\epsilon}^{p}(\boldsymbol{d}_{n_3}) \right)\\ \nonumber
	+ &\frac{f(\boldsymbol{x}_{s_2}) \pm f(\boldsymbol{x}_{p})}{\epsilon(\boldsymbol{x}_{p})}  \left( \eta_{\epsilon}^{p}(\boldsymbol{d}_{s_2}) +  \eta_{\epsilon}^{p}(\boldsymbol{d}_{n_4}) + \eta_{\epsilon}^{p}(\boldsymbol{d}_{n_5}) + \eta_{\epsilon}^{p}(\boldsymbol{d}_{n_6}) \right) \\
	+ &\frac{f(\boldsymbol{x}_{p}) \pm f(\boldsymbol{x}_{p})}{\epsilon(\boldsymbol{x}_{p})}  \left( \eta_{\epsilon}^{p}(\boldsymbol{d}_{n_7}) +  \eta_{\epsilon}^{p}(\boldsymbol{d}_{n_8}) + \eta_{\epsilon}^{p}(\boldsymbol{d}_{n_9}) + \eta_{\epsilon}^{p}(\boldsymbol{d}_{n_{10}}) \right).
\end{align}
Since the embedding-space kernels are evaluated at concrete distances, the $\eta_{\epsilon}^{p}(\boldsymbol{d}_q)$ are just scalar numbers.
All kernel values that share the same pre-factor $\frac{f(\boldsymbol{x}_q) \pm f(\boldsymbol{x}_p)}{\epsilon(\boldsymbol{x}_p)}$ are thus summed to the surface kernels $\eta_{\mathcal{S}}(\boldsymbol{x}_p,\boldsymbol{x}_q)$, $q \in \mathcal{N}_{\mathcal{S}}$, obtaining:
\begin{align}
	\nonumber
	\boldsymbol{Q}_\mathcal{S}^{\boldsymbol{\alpha}} f(\boldsymbol{x}_{p}) =
	&\frac{f(\boldsymbol{x}_{s_1}) \pm f(\boldsymbol{x}_{p})}{\epsilon(\boldsymbol{x}_{p})} \; \eta_{\mathcal{S}}(\boldsymbol{x}_p,\boldsymbol{x}_{s_1}) +
	\frac{f(\boldsymbol{x}_{s_2}) \pm f(\boldsymbol{x}_{p})}{\epsilon(\boldsymbol{x}_{p})} \; \eta_{\mathcal{S}}(\boldsymbol{x}_{p},\boldsymbol{x}_{s_2}) \\
	+ &\frac{f(\boldsymbol{x}_{p}) \pm f(\boldsymbol{x}_{p})}{\epsilon(\boldsymbol{x}_{p})} \; \eta_{\mathcal{S}}(\boldsymbol{x}_{p},\boldsymbol{x}_{p}).
\end{align}
For even differential operators, i.e., derivatives with even $|\boldsymbol{\alpha}|$, the third term vanishes identically and can be skipped in the calculations. But this is not the case for odd-order derivatives.
With this rearrangement, each evaluation of the operator at a point $p$ only requires three kernel evaluations instead of the 12 that would be required in the embedding case. In addition, the normally extended points never need to be allocated and stored, as all kernel computations can happen on the fly. Algorithm~\ref{alg:surface_kernels} details the procedure for Surface DC-PSE operator construction. For the example from Fig.~\ref{fig:sdcpse}b, this results in the surface DC-PSE kernel values:
$$
\mathcal{K}_{\mathcal{S}} = [
\eta_{\mathcal{S}}(\boldsymbol{x}_p,\boldsymbol{x}_{s_1}),\,
\eta_{\mathcal{S}}(\boldsymbol{x}_p,\boldsymbol{x}_{s_2}),\,
\eta_{\mathcal{S}}(\boldsymbol{x}_p,\boldsymbol{x}_p)
   ],
$$
which can directly be used in Eq.~(\ref{eq:sdcpseEval}) to evaluate surface differential operators.

\begin{algorithm}
	%	\SetAlgoLined
	\textbf{Input:}
	\begin{enumerate}
      \item List of neighbor distances $\mathcal{N}_{\text{dist}}$ for each pair $\boldsymbol{x}_{p,q}$, with $q \in \mathcal{N}(\boldsymbol{x}_p)$, as constructed by Algorithm \ref{alg:embedding_neighborhood}
	\end{enumerate}
	\textbf{Output:} Surface kernel values for each pair $\boldsymbol{x}_{p,s}$, with $s \in \mathcal{N}_\mathcal{S}$: $\mathcal{K}_{\mathcal{S}}$
	\\
	\begin{algorithmic}[1]
		\State  $\mathcal{K}_{\mathcal{S}} = zeros(size(\mathcal{N}_{\mathcal{S}}) + 1)$
      \State $\eta = \text{DCPSE}(p)$ \Comment{Determine embedding-space DC-PSE kernel at $p$}
      \ForAll{$(i,j) \in \mathcal{N}_{\text{dist}}$}
         \State $\mathcal{K}_{\text{emb}}[i][j] = \eta(\mathcal{N}_{\text{dist}}[i][j])$
      \EndFor
		\For{$i \in [1,size(\mathcal{N}_{\mathcal{S}})+1]$}
      \ForAll{$j \in \mathcal{K}_{\text{emb}}[i]$}
      \State $\mathcal{K}_{\mathcal{S}}[i] \mathrel{+}= \mathcal{K}_{\text{emb}}[i][j]$
			\EndFor
		\EndFor
		\Statex
	\end{algorithmic}
	\caption{Surface DC-PSE: construction of surface kernel at a point $p$.}
	\label{alg:surface_kernels}
\end{algorithm}

\section{Results} \label{sec:Result}
We validate and benchmark the Surface DC-PSE method. First, we verify its convergence in test cases with known analytical solution. Then, we show applications to cases with more general surfaces where no analytical solution is available. Finally, we compare with surface finite-element methods (SFEM) in a validation study.

In all cases, orthogonal extension is exact, copying surface points along the known normals.
Therefore, $\nabla F \cdot \boldsymbol{n}$ in Eq.~\ref{eq:surfaceIdentity} is always zero by construction. Numerically evaluating this term requires approximating the gradient $\nabla F$, which we confirmed to converge with the order of accuracy of the discretization scheme used.

\subsection{Laplace-Beltrami operator on a circle and a sphere}
\begin{figure}[t]
	\setlength{\tabcolsep}{-5pt}
	\begin{tabular}{cc}
		\subfloat[]{\includegraphics[scale=0.175]{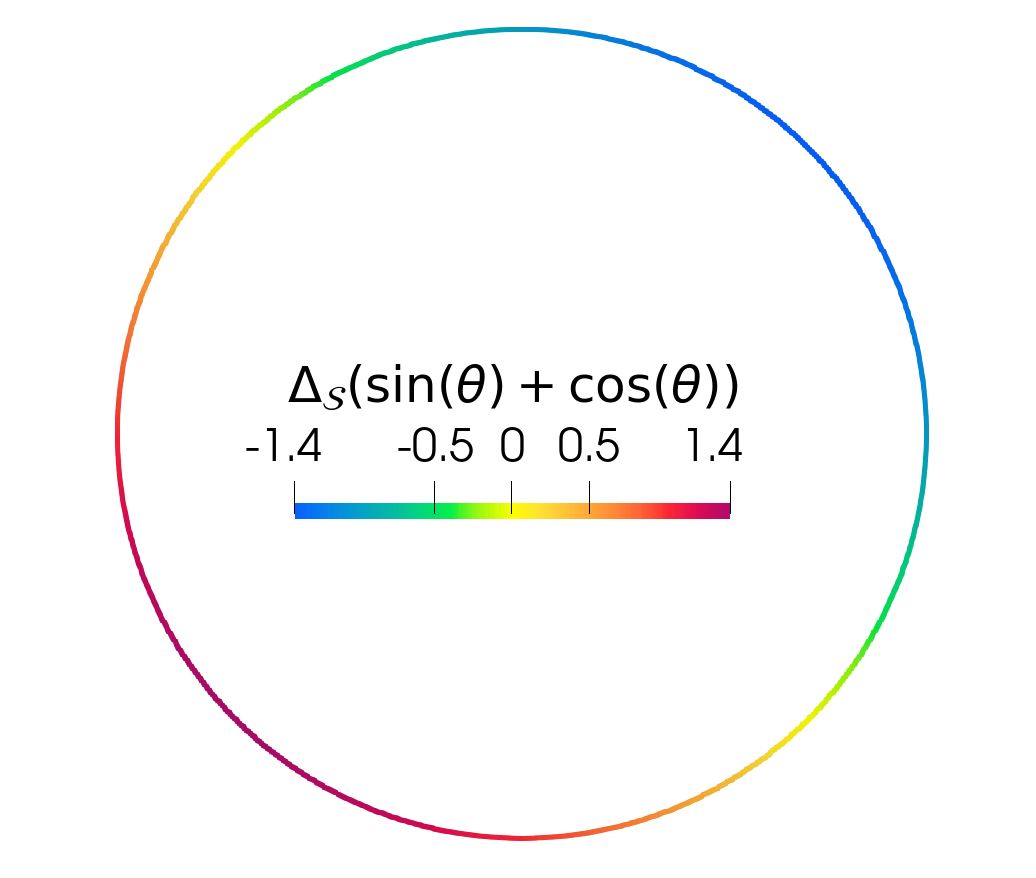}} & \subfloat[]{\includegraphics[scale=0.175]{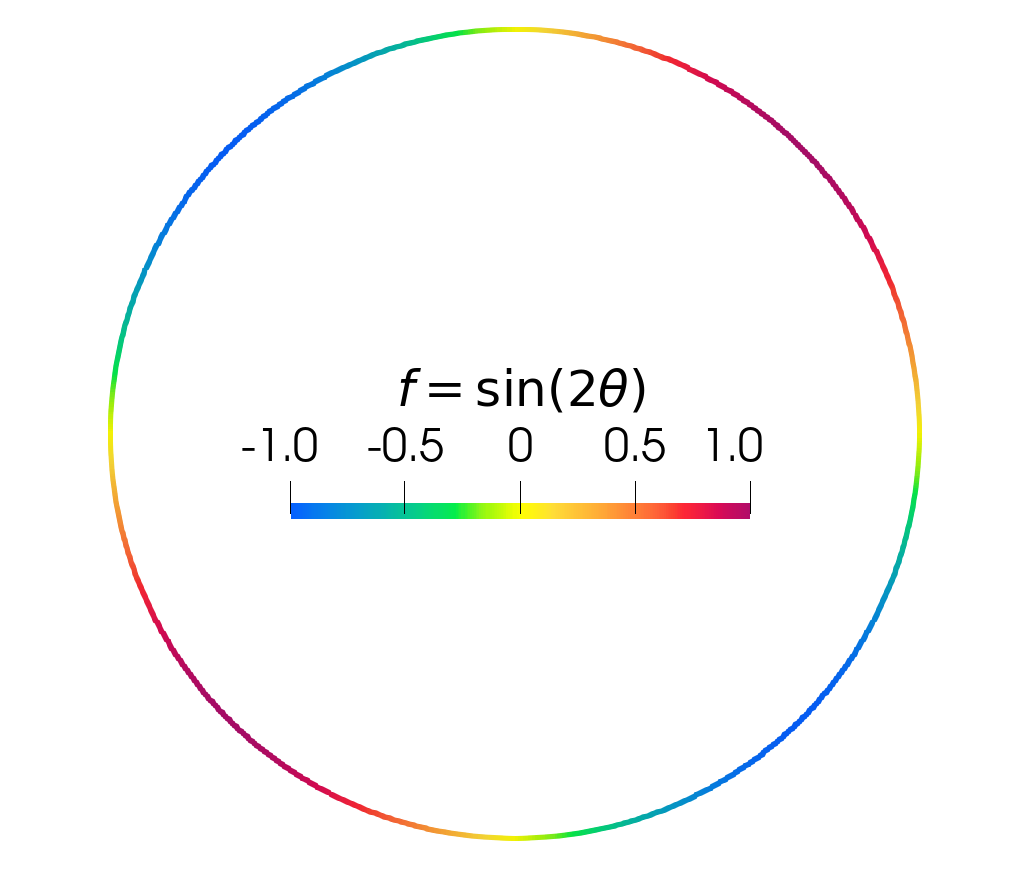}} \\
		\subfloat[]{\includegraphics[scale=0.34]{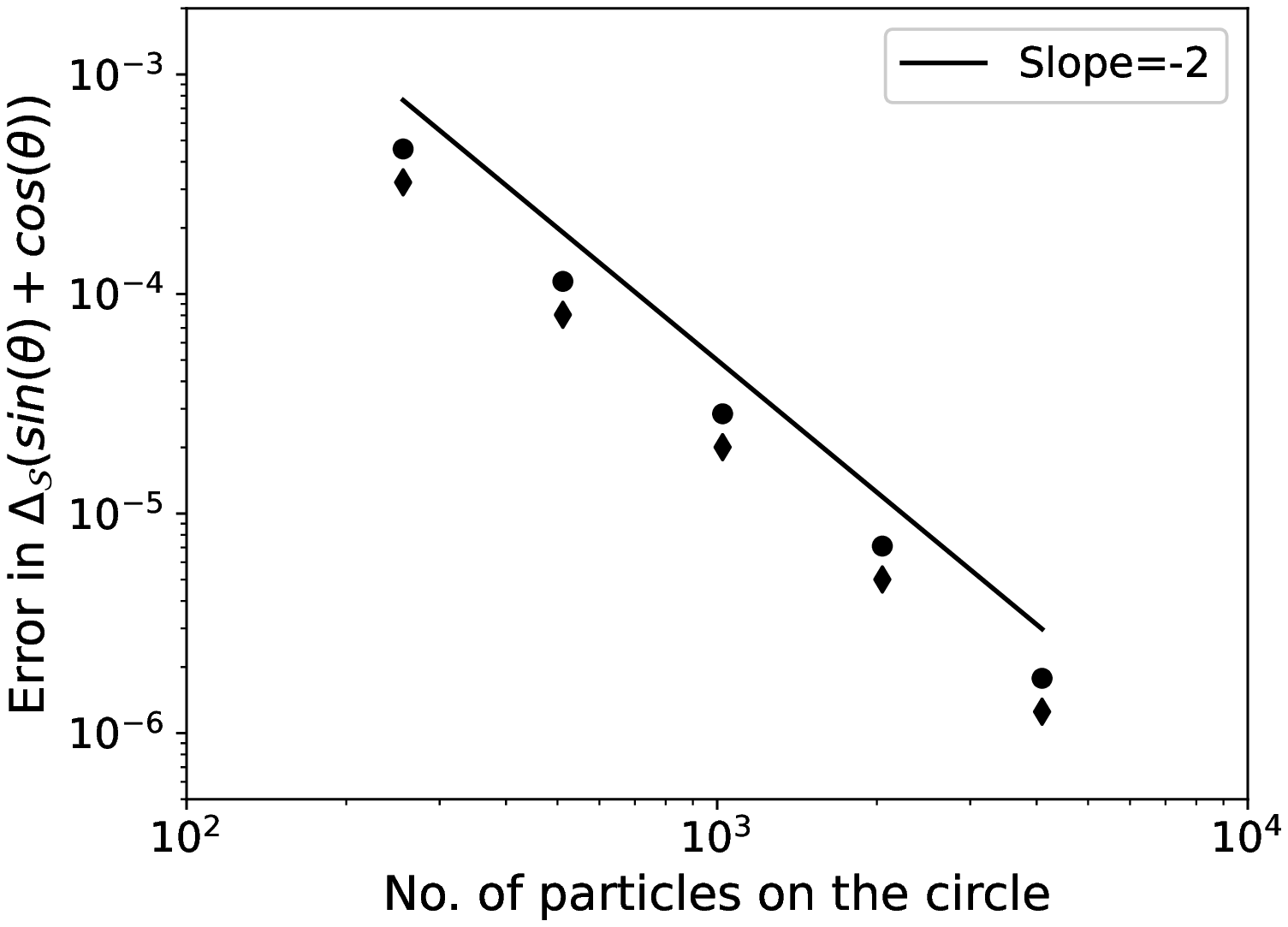}}	 & \subfloat[]{\includegraphics[scale=0.34]{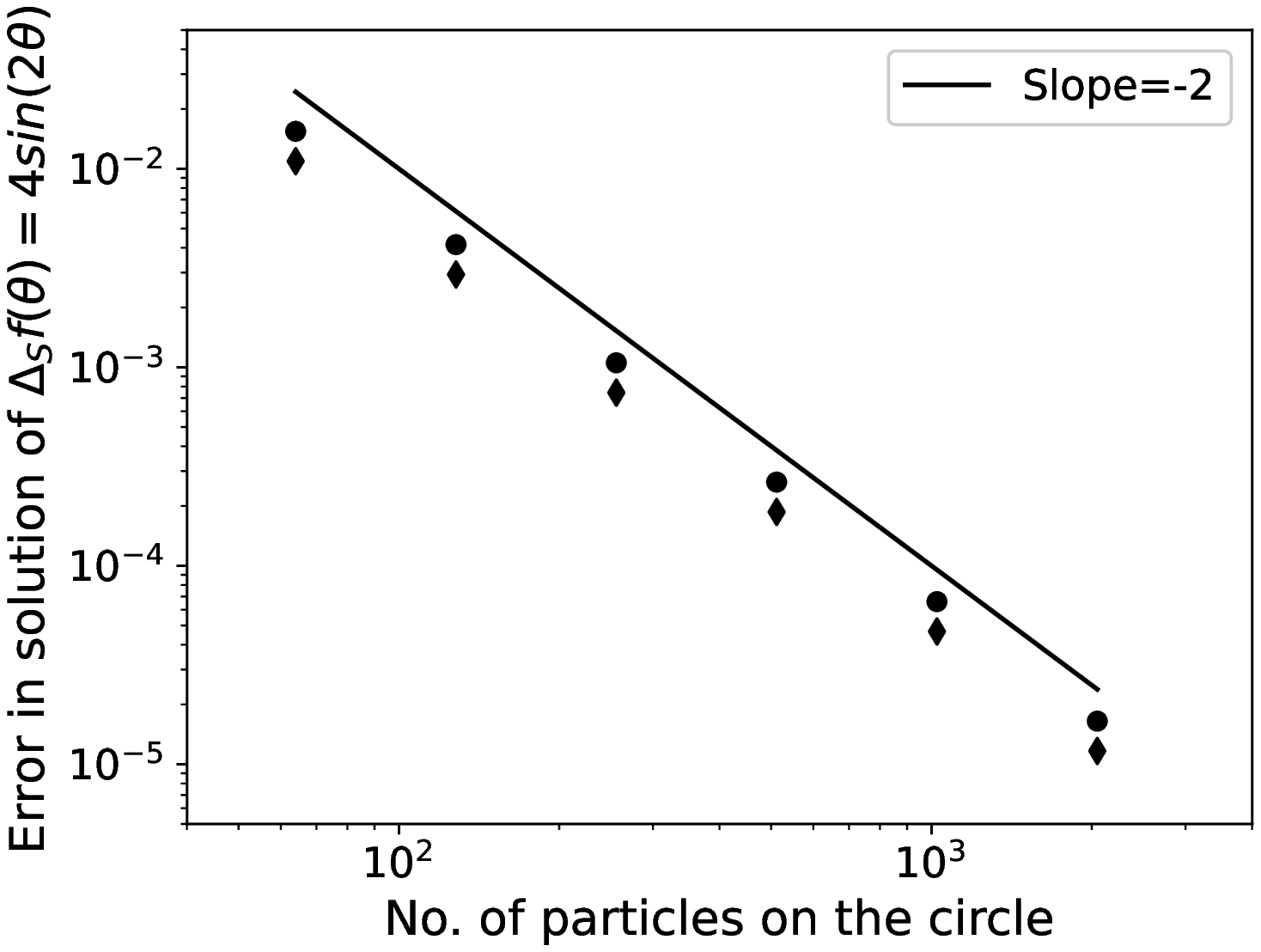}}
	\end{tabular}
	\caption{\textbf{Visualization and convergence of the Laplace-Beltrami operator on the unit circle.
			(a)} Visualization of $\Lap_{\mathcal{S}}(\sin(\theta)+\cos(\theta))$ computed using second-order accurate Surface DC-PSE operators.
		\textbf{(b)} Visualization of the solution of the Poisson equation $\Lap_{\mathcal{S}}f=4\sin(2\theta)$ solved using second-order accurate Surface DC-PSE  operators.
		\textbf{(c)} Convergence plot of the Laplace-Beltrami operator in (a). $L_\infty(\bullet)$ and $L_2(
		\blacklozenge)$  norms of the absolute errors are computed against the analytical solution in Eq.~(\ref{eq:circLBSol}) for increasing numbers of points on the circle.
		\textbf{(d)} Convergence plot of the Poisson equation solution in (b). $L_\infty(\bullet)$ and $L_2 (\blacklozenge)$  norms of the absolute errors are computed against the analytical solution in Eq.~(\ref{eq:circPoisSol}) for increasing numbers of points on the circle.}
\label{fig:LBCirc}
\end{figure}

We start by verifying convergence for the Laplace-Beltrami operator on the unit circle  $S^1$.
The collocation points are distributed regularly using equi-angular spacing. We use a normal spacing of $\delta n = 3/(N_p-1)$ to compute the surface operators in Eq.~(\ref{eq:sdcpseEval}) in a narrow band with $N_n = 4$ layers on each side of the surface.
$N_p$ is the number of surface points $\boldsymbol{x}_s$. We choose $r_c = 4.1\delta n$ as the operator support and $r=2$ as the desired order of convergence. The Laplace-Beltrami operator is characterized by the multi-index $\boldsymbol{\alpha}=(2,0)+(0,2)$. Note that this multi-index is 2-dimensional, despite the circle being one-dimensional, since the operators are constructed in the embedding space, but evaluated intrinsically.

We test the numerical approximation of the surface operator on the function
\begin{equation}
f(\theta)=\sin(\theta)+\cos(\theta)
\end{equation}
in polar coordinates.
The error is computed against the analytical solution
\begin{equation}
\mathrm{\Delta}_{\mathcal{S}} f(\theta)= \nabla_{\!\!\mathcal{S}}\cdot (\nabla_{\!\!\mathcal{S}}f(\theta)) =\nabla_{\!\theta}^2 f(\theta)=-(\sin(\theta)+\cos(\theta)).
\label{eq:circLBSol}
\end{equation}
The visualization of the numerical solution and the convergence plot of the absolute errors are shown in Fig.~\ref{fig:LBCirc}a,c. We observe second-order convergence to the analytical solution, as expected for $r=2$.

\begin{figure}[h!]
\setlength{\tabcolsep}{0pt}
\centering
\begin{tabular}{cc}
\centering
\subfloat[]{\includegraphics[width=0.48\textwidth,scale=0.15,trim=0 20 0 0, clip]{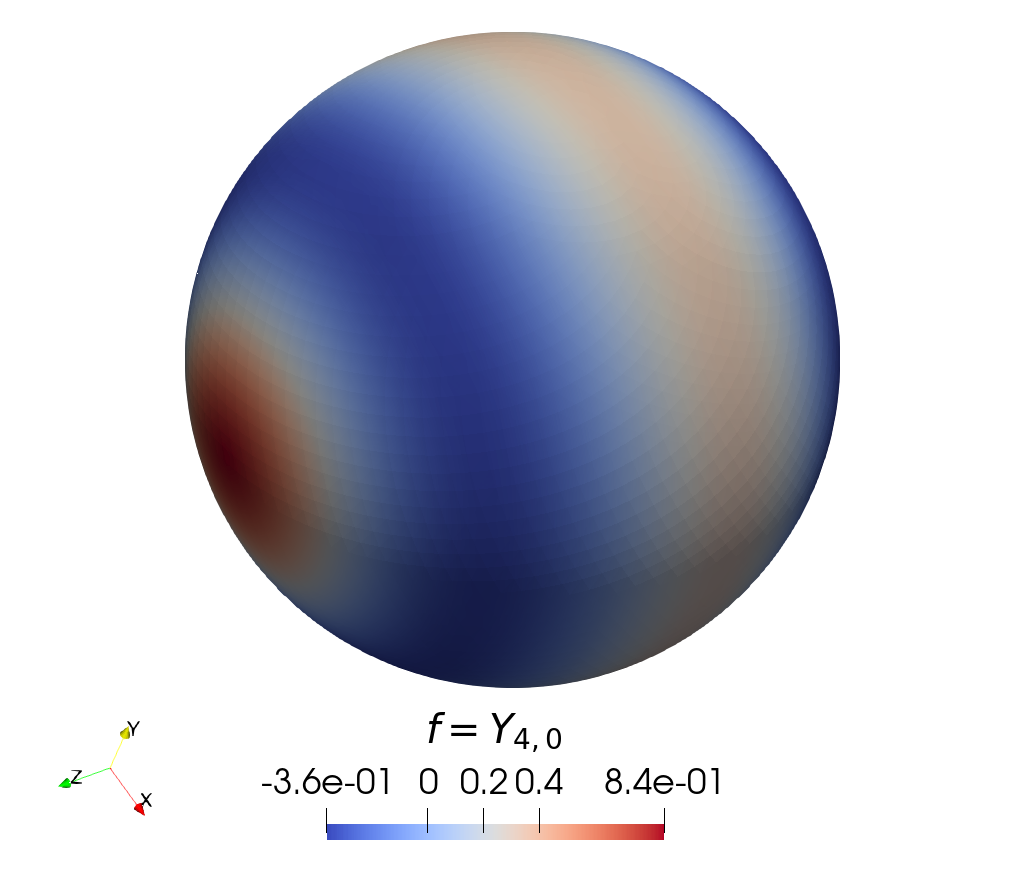}} & 	\subfloat[]{\includegraphics[width=0.45\textwidth,scale=0.4, trim= 0 -10 -20 0]{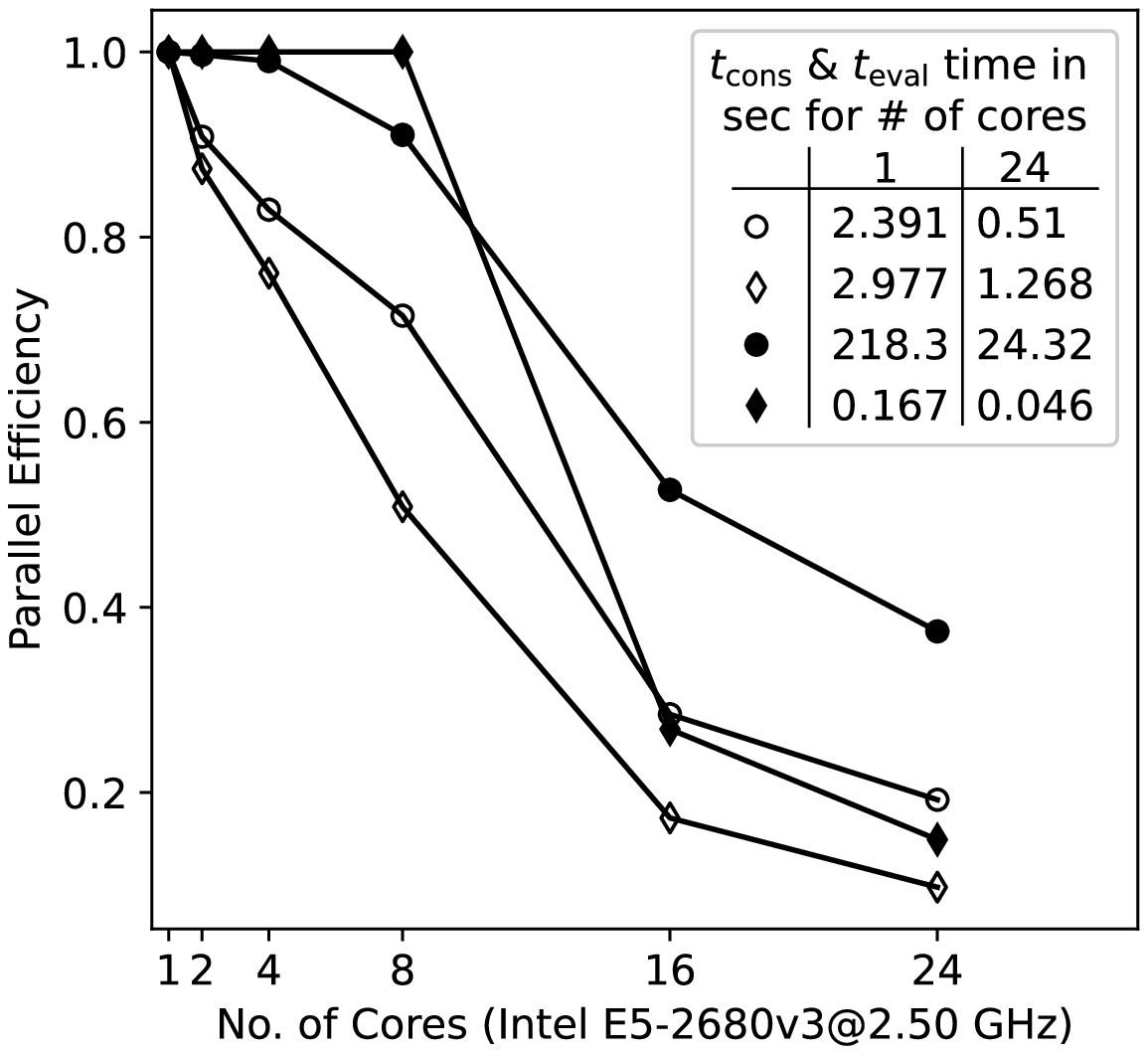}} \\
\subfloat[]{\includegraphics[width=0.465\textwidth,scale=0.42]{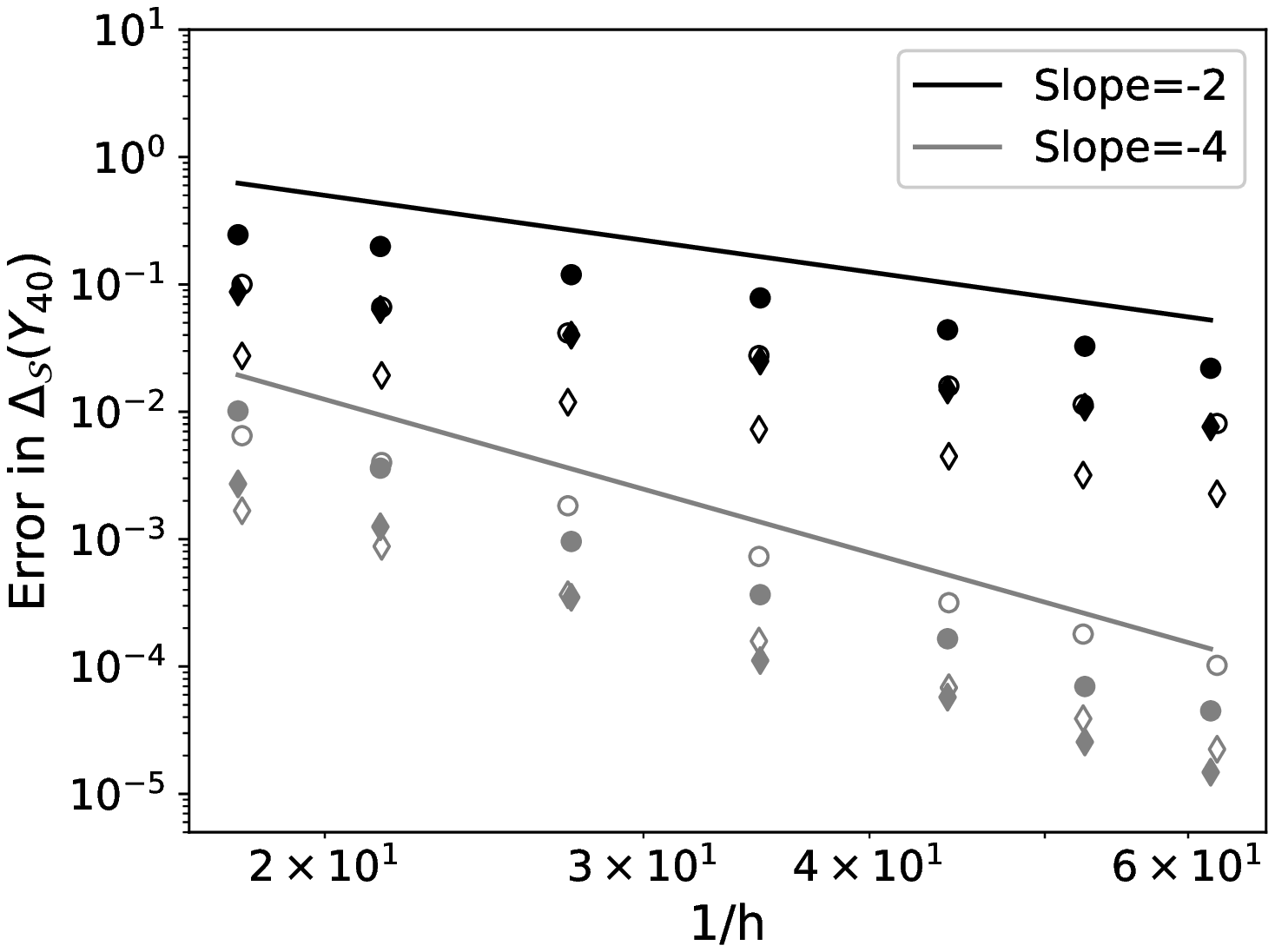}} & 			\subfloat[]{\includegraphics[width=0.49\textwidth,scale=0.42]{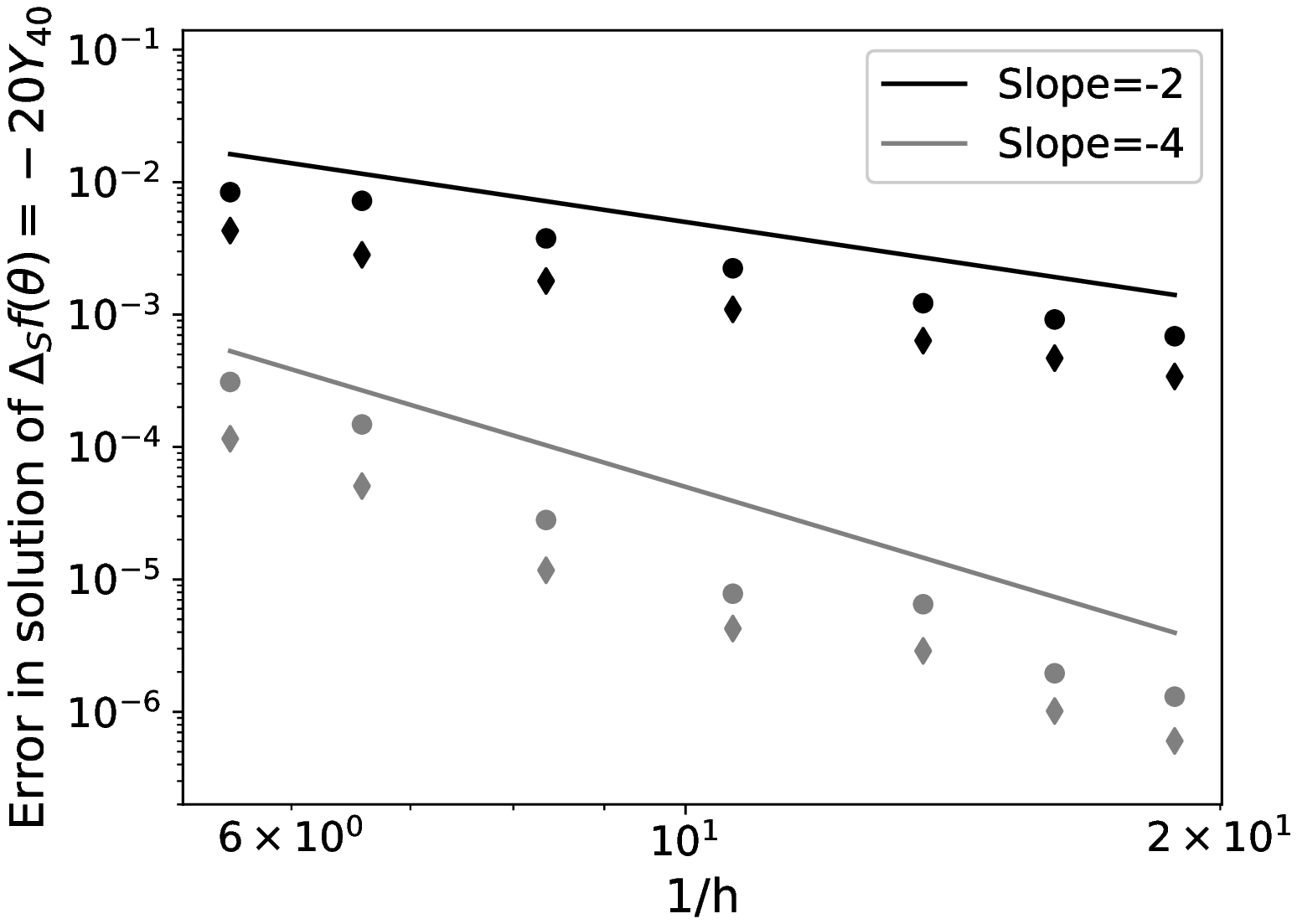}}
\end{tabular}
\caption{\textbf{Visualization and convergence of the Laplace-Beltrami operator on the unit sphere.
(a)} Visualization of the spherical harmonic function $Y_{4,0}$.
\textbf{(b)} Results of a strong scalability test for the computation of $\Lap_{\mathcal{S}}Y_{4,0}$ for average spacing $h=0.05$ using Surface DC-PSE operators $(\bullet,\blacklozenge)$ and the Closest Point (CP) method ($\circ,\lozenge$)~\cite{ruuth_simple_2008} on increasing numbers of CPU cores. We separately show the parallel efficiency for the construction of
Surface DC-PSE operators $(\bullet)$ and of the CP representation $(\circ)$, and for their evaluation to approximate
$\Lap_{\mathcal{S}}Y_{4,0}$ across the entire domain using Surface DC-PSE $(\blacklozenge)$ and the CP method $(\lozenge)$. The absolute wall-clock times for one time step are shown in the inset table for 1 and 24 cores.
\textbf{(c)} Convergence plot of $\Lap_{\mathcal{S}}Y_{4,0}$ for Surface DC-PSE ($L_\infty (\bullet)$, $L_2 (\blacklozenge)$) and the CP method ($L_\infty (\circ)$, $L_2 (\lozenge)$) using second-order (black) and fourth-order (gray) approximations. Norms of the absolute errors are computed against the analytical solution in Eq.~(\ref{eq:sphLBSol}) for increasing numbers of points (decreasing average spacing $h$).
\textbf{(d)} Convergence plot for the solution of the Poisson equation $\Lap_{\mathcal{S}}f=-20Y_{4,0}$ across the entire domain using second-order (black) and fourth-order (gray) operators.
$L_\infty (\bullet)$  and $L_2 (\blacklozenge)$  norms of the absolute errors are computed against the analytical solution in Eq.~(\ref{eq:sphPoisSol}) for increasing numbers of surface points.}
\label{fig:SphPois}
\end{figure}

We further test the method on the unit sphere $S^2$. The collocation points are distributed using the Fibonacci sphere technique \cite{Gonzlez2009}. We use a normal spacing of $\delta n = 0.8/(\sqrt[3]{N_p}-1)$ to determine the surface operators in Eq.~(\ref{eq:sdcpseEval}) in a narrow band with $N_n=2$ layers on each side of the surface.
$N_p$ is the number of points on the sphere. We choose $r_c = 2.9\delta n$ as the operator support and $r=2$ and $r=4$ as the desired orders of convergence.
The Laplace-Beltrami operator is characterized by the three-dimensional multi-index $\boldsymbol{\alpha}=(2,0,0)+(0,2,0)+(0,0,2)$.
We test the numerical approximation of the surface operator on the scalar spherical harmonic function
\begin{equation}
f(\theta,\phi)=Y_{lm}
\end{equation}
in spherical coordinates for the mode $l=4$, $m=0$ (Fig. \ref{fig:SphPois}a).
The error is computed against the analytical solution
\begin{equation}
\mathrm{\Delta}_{\mathcal{S}} f= \nabla_{\!\!\mathcal{S}}\cdot (\nabla_{\!\!\mathcal{S}}f(\theta,\phi)) =-l(l+1)Y_{lm}
\label{eq:sphLBSol}
\end{equation}
and is plotted in Fig.~\ref{fig:SphPois}c.

We also use this test case to benchmark against the Closest Point (CP) method \cite{ruuth_simple_2008} with $L_2$ and $L_\infty$ errors plotted in Fig.~\ref{fig:SphPois}c. While Surface DC-PSE is less accurate than the CP method for second-order operators, it outperforms for fourth-order operators, which is likely due to the better condition numbers of the small linear systems to be solved for each point.

Finally, we perform a strong scaling benchmark of the computation time with increasing numbers of CPU cores with both codes implemented in the parallel computing library OpenFPM~ \cite{incardona_openfpm_2019,singh_c_2021} in C++ and run on the same hardware.
We plot the parallel efficiency (i.e., the speed-up divided by the number of cores) in Fig.~\ref{fig:SphPois}b for both the construction and the evaluation of the operators of both methods. We further report the absolute wall-clock times on one and 24 cores in the inset table. These times show that evaluating the Surface DC-PSE operators over all points in the domain is about one order of magnitude faster than evaluating the CP transform~\cite{ruuth_simple_2008}. In addition, Surface DC-PSE scales better with increasing numbers of parallel CPU cores. This is due to the simpler kernel evaluation of DC-PSE\footnote{Up to 8 cores, the measured parallel efficiency of Surface DC-PSE was larger than 1.0, due to cache effects when storing the kernel coefficients.}. Eventually, the efficiency for both methods drops, as is expected for strong scaling with constant communication overhead.
The time required to construct the Surface DC-PSE operators with $N_n=2$, however, is about two orders of magnitude larger than that for constructing the CP representation in a narrow band of radius 5 as required by the regression support. However, it still scales better with increasing CPU core count.
For Eulerian simulations, where collocation points do not move, the kernels have to be determined only once at the beginning of a simulation, or they can be loaded from files for standard point distributions, potentially leading to an overall lower computational cost of Surface DC-PSE.

\subsection{Poisson equation on a circle and a sphere}

Surface DC-PSE operators can also be used for implicit equations by solving a linear system of equations with a system matrix constructed using the Surface DC-PSE operators. We test this by solving the Poisson equation on the unit circle $S^1$:
\begin{equation}
\Lap _{\mathcal{S}} f = 4 \sin(2\theta)\qquad\theta\in \mathrm{\Omega}=S^1\backslash (1,0)
\end{equation}
with Dirichlet boundary condition at one point $(1,0)$ conforming to the analytical solution
\begin{equation}
f (\theta) = \sin(2\theta)\qquad \theta\in \mathrm{\Omega}\cup (1,0).
\label{eq:circPoisSol}
\end{equation}
We use the same Surface DC-PSE operators as in the previous subsection to construct the system of equations, which is then solved using the KSPGMRES solver from PETSc~\cite{balay_efficient_1997}. Fig.~\ref{fig:LBCirc}b,d show the solution $f$ and the convergence plot of the absolute error with respect to the analytical solution.

Next, we test the method in three dimensions by solving the Poisson equation on the sphere $S^2$:
\begin{equation}
\Lap _{\mathcal{S}} f = -20 Y_{4,0}
\label{eq:sphPois}
\end{equation}
with Dirichlet boundary condition along the great circle parallel to the $y-z$ plane conforming to the analytical solution
\begin{equation}
f = Y_{4 0}.
\label{eq:sphPoisSol}
\end{equation}
We solve the resulting linear system with KSPGMRES from PETSc~\cite{balay_efficient_1997} without preconditioning. The convergence plots for orders $r=2$ and $r=4$ are shown in Fig.~\ref{fig:SphPois}d. The collocation points on the sphere were constructed using the Fibonacci sphere technique \cite{Gonzlez2009}. While this generates pseudo-regular point distributions on the sphere, their average spacing does fluctuate a bit, explaining the slight waviness of the curve especially for the fourth-order operators.

\subsection{Mean and Gauss curvature computation}

\begin{figure}[h!]
\setlength{\tabcolsep}{0pt}
\centering
\begin{tabular}{cc}
\subfloat[]{\includegraphics[width=0.49\textwidth,scale=0.175]{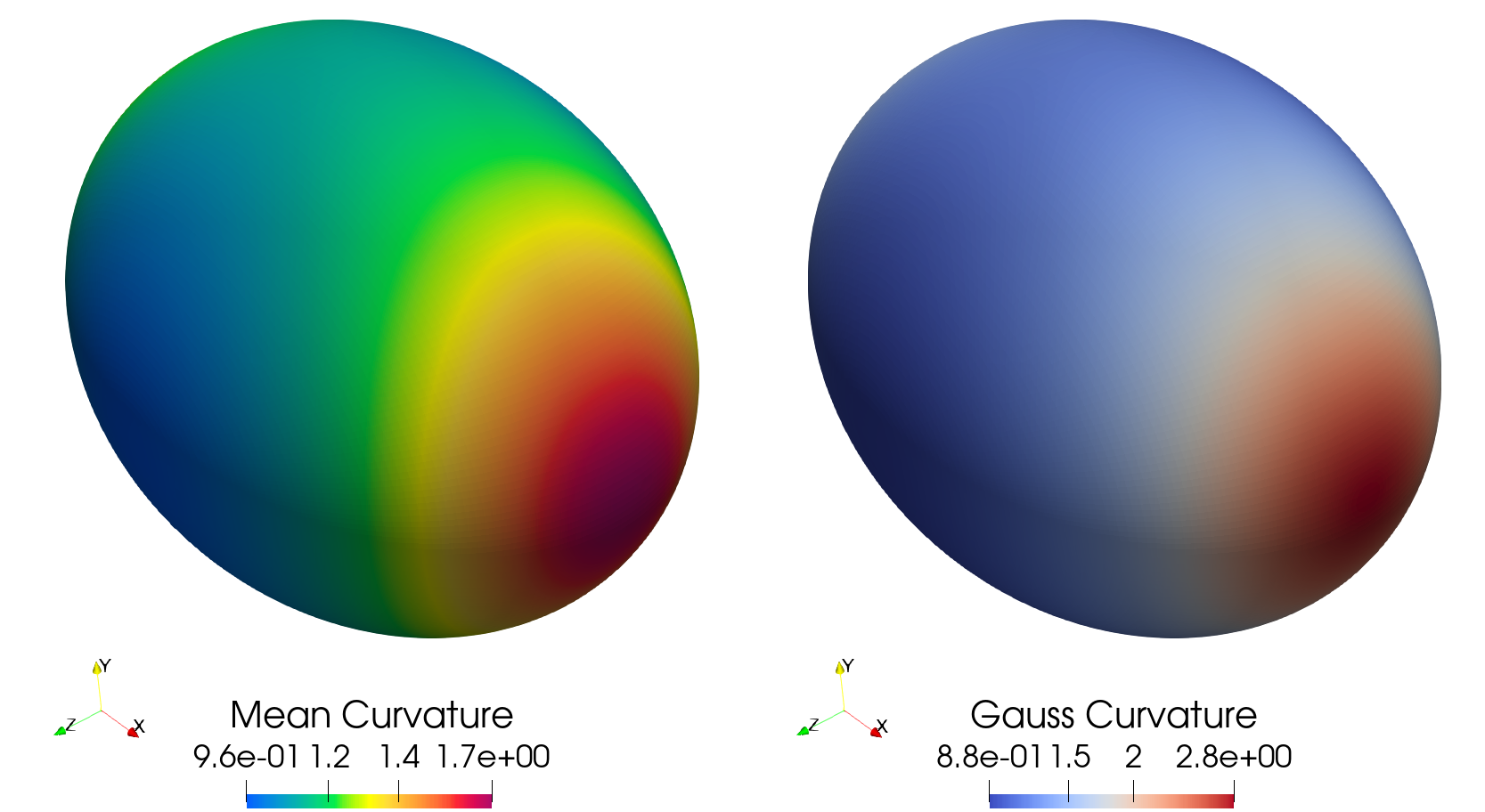}} & 			\subfloat[]{\includegraphics[width=0.49\textwidth,scale=0.175]{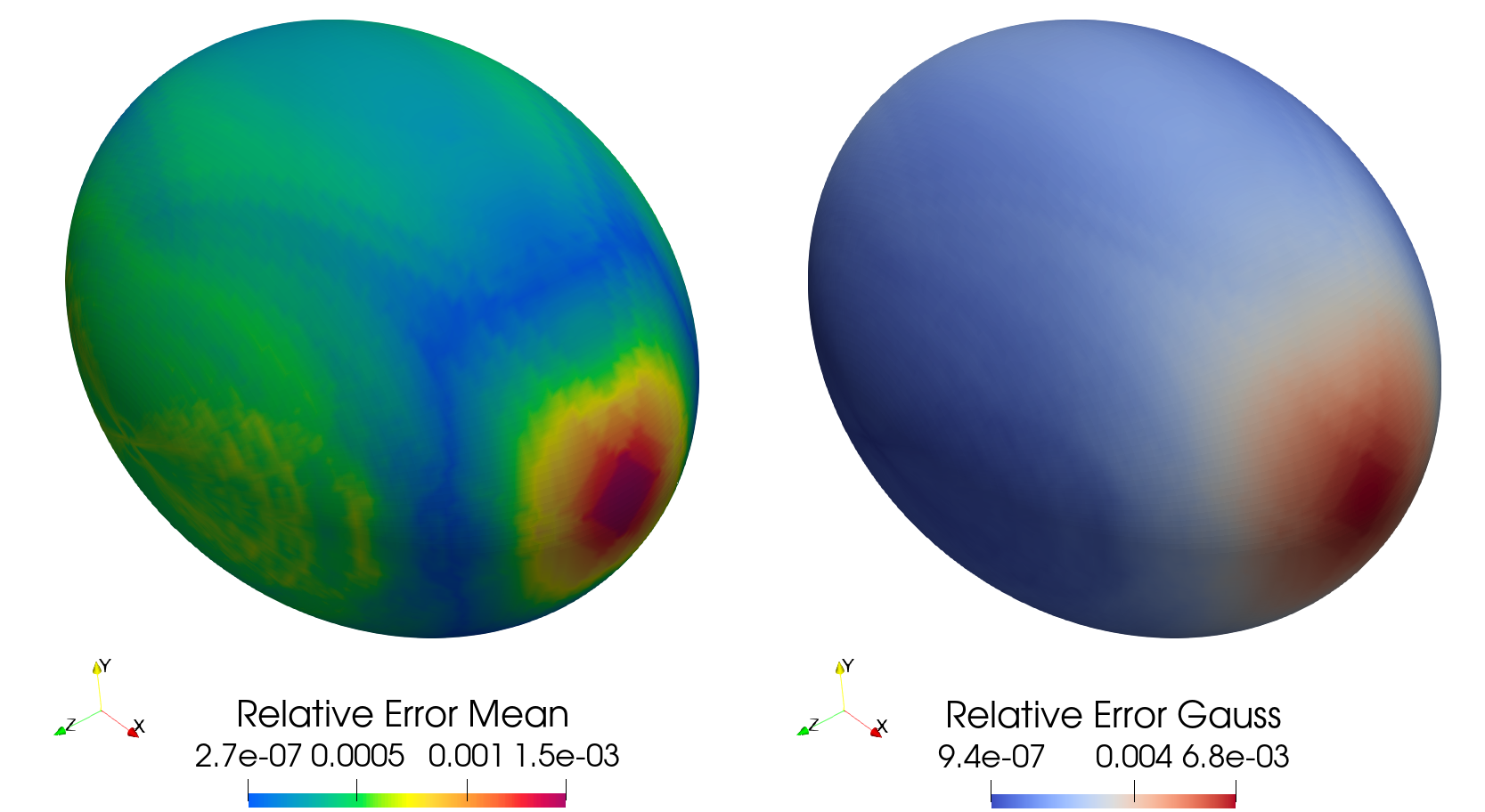}} \\
\subfloat[]{\includegraphics[width=0.49\textwidth,scale=0.3]{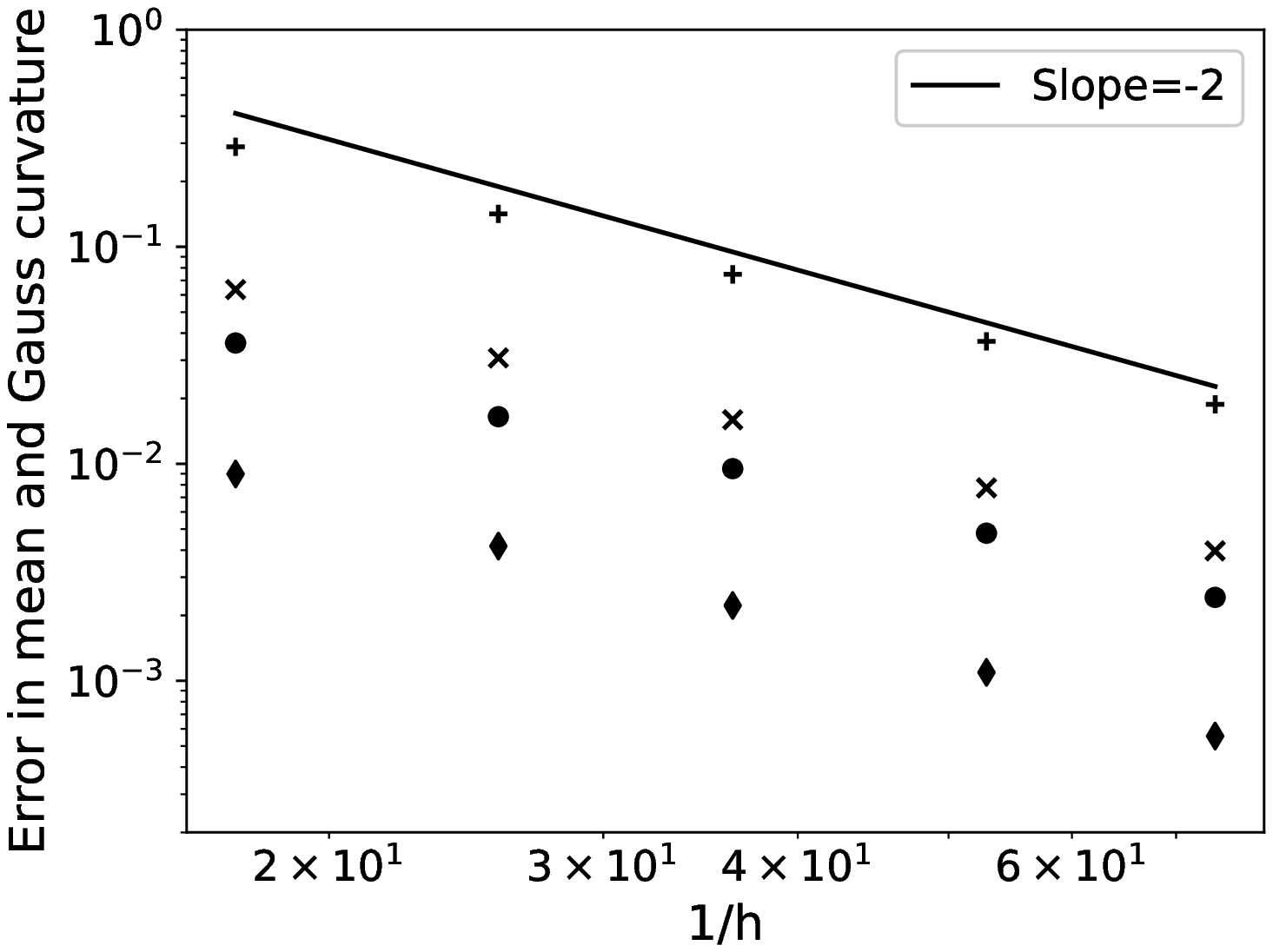}}	 & \subfloat[]{\includegraphics[width=0.49\textwidth,scale=0.09]{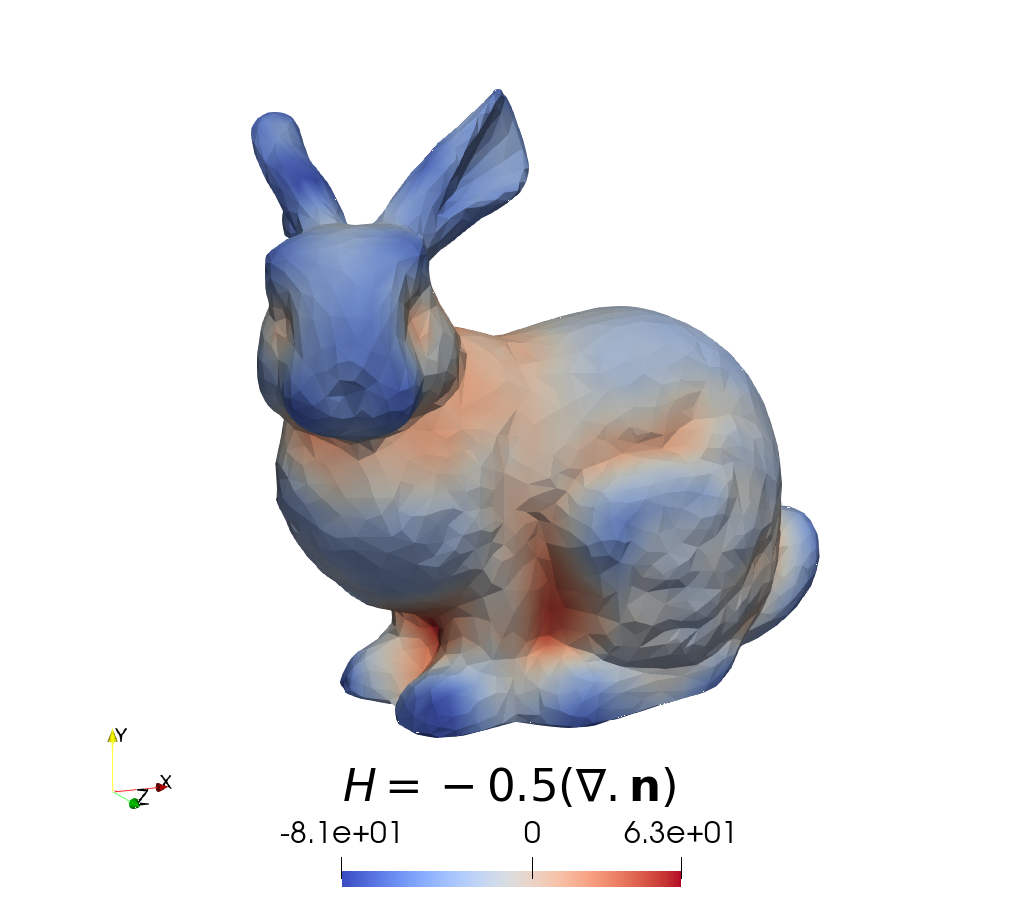}}
\end{tabular}
\caption{\textbf{
Gauss and mean curvature computation using Surface DC-PSE.
(a)} Visualization of the mean and Gauss curvatures of an ellipsoid, numerically computed as the trace of the embedded shape tensor $-0.5\nabla_{\!\mathcal{S}}\cdot\boldsymbol{n}=-0.5\mathrm{Tr}(\nabla_{\!\mathcal{S}}\boldsymbol{n})$ of the surface normal vector field $\boldsymbol{n}$ and the product of the non-zero eigenvalues of $\nabla_{\!\mathcal{S}}\boldsymbol{n}$, respectively, using second-order accurate Surface DC-PSE operators. \textbf{(b) } Visualization of the relative errors in the computed curvatures from (a) on 32258 surface points in comparison with the analytical solution in Eq.~(\ref{eq:EllipErr}). They range from about $10^{-7}$ to $10^{-3}$. \textbf{(c)} Convergence plot of the mean and Gauss curvature computations. $L_\infty(\bullet,+)$  and $L_2(\blacklozenge,\times)$ absolute errors for mean and Gauss curvatures, respectively, are computed against the analytical solutions in Eq.~(\ref{eq:EllipErr}) for decreasing average collocation point spacing $h$. \textbf{(d)} Visualization of the mean curvature computed on the Stanford bunny with 2960 points using second-order accurate Surface DC-PSE operators.}
\label{fig:ellipCurv}
\end{figure}

We further verify Surface DC-PSE by computing the mean curvature $H$ and Gauss curvature $K$ of an ellipsoid
\begin{equation}
\frac{x^2}{a^2}+\frac{y^2}{b^2}+\frac{z^2}{c^2}=1
\end{equation}
with $a=1$, $b=0.8$, $c=0.75$ and parameterization $(u,v)$
\begin{equation}
x =a \cos u \sin v, \qquad	y =b \sin u \sin v, \qquad	z =c \cos v.
\end{equation}
We compute both curvatures from the embedded shape tensor $\nabla_{\!\mathcal{S}}\boldsymbol{n}$, i.e.,~the extension of the intrinsic shape operator to the embedding space. We numerically compute this tensor in Cartesian coordinates as the $3\times 3$ matrix
$$
 \nabla_{\!\mathcal{S}}\boldsymbol{n} =
 \left( \begin{array}{c}
 \nabla_{\!\mathcal{S}}n_1 \\
 \nabla_{\!\mathcal{S}}n_2 \\
 \nabla_{\!\mathcal{S}}n_3
\end{array}\right),
$$
where $\boldsymbol{n} = (n_1, n_2, n_3)$ are the components of the analytically given normal vectors at the collocation point.
All intrinsic derivatives over the scalar fields $n_1, n_2, n_3$ are approximated by Surface DC-PSE operators. Mean curvature $H$ is then computed as the trace of the embedded shape tensor and Gauss curvature $K$ as the product of its non-zero eigenvalues (principal curvatures). These numerical computations are verified against the analytical solutions
\begin{align}
    H &=\frac{abc[3(a^2+b^2)+2c^2+(a^2+b^2-2c^2)\cos(2v)-2(a^2-b^2)\cos(2u)\sin^2v]}{8[a^2b^2\cos^2v+c^2(b^2\cos^2u+a^2\sin^2u)\sin^2v]^{3/2}}\notag
\\
   K &=\frac{a^2 b^2 c^2}{\left[a^2 b^2 \cos ^2 v+c^2\left(b^2 \cos ^2 u+a^2 \sin ^2 u\right) \sin ^2 v\right]^2}.
\label{eq:EllipErr}
\end{align}
We approximate the embedded shape tensor $\nabla_{\!\mathcal{S}} \boldsymbol{n}$ using Surface DC-PSE with $\delta n =3.0/(N_p-1)$, $r_c=2.9 \delta n$, $N_n=2$, and $r=2$. The results and the convergence plot of the absolute errors are shown in Fig.~\ref{fig:ellipCurv}a,c. As specified by $r$, we observe second-order convergence to the analytical solution when decreasing the average spacing $h$ between the points. The relative errors are visualized in Fig.~\ref{fig:ellipCurv}b. They concentrate around extremal points of the curvature, as expected.

Finally, we apply the same mean-curvature computation to an arbitrary surface with no analytical solution, the Stanford bunny from the Stanford Computer Graphics Laboratory.
We use the down-sampled version of the original data set with 2960 points on the surface, obtained from \url{https://www.stlfinder.com/model/stanford-bunny-S4kAUsKI/3091553}.
The result is visualized in Fig.~\ref{fig:ellipCurv}d, showing that the Surface DC-PSE qualitatively works also for non-algebraic surfaces.

\subsection{Comparison with Surface Finite Element Methods}\label{sec:fem}

\begin{figure}[h!]
\setlength{\tabcolsep}{0pt}
\centering
\begin{tabular}{cc}
	\subfloat[]{\includegraphics[width=0.49\textwidth,scale=0.19]{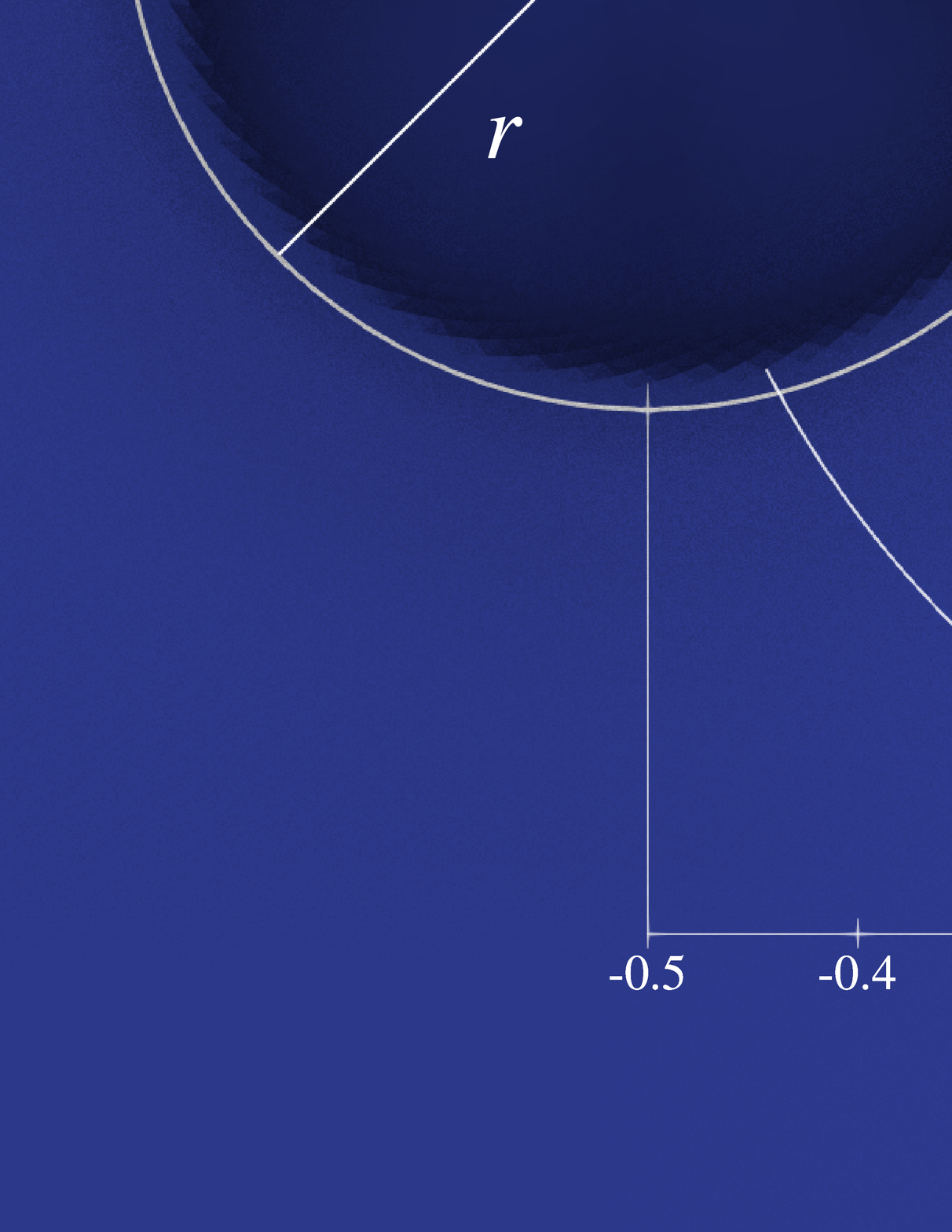}} &
	 \quad
	 \subfloat[]{\includegraphics[width=0.49\textwidth,scale=0.19]{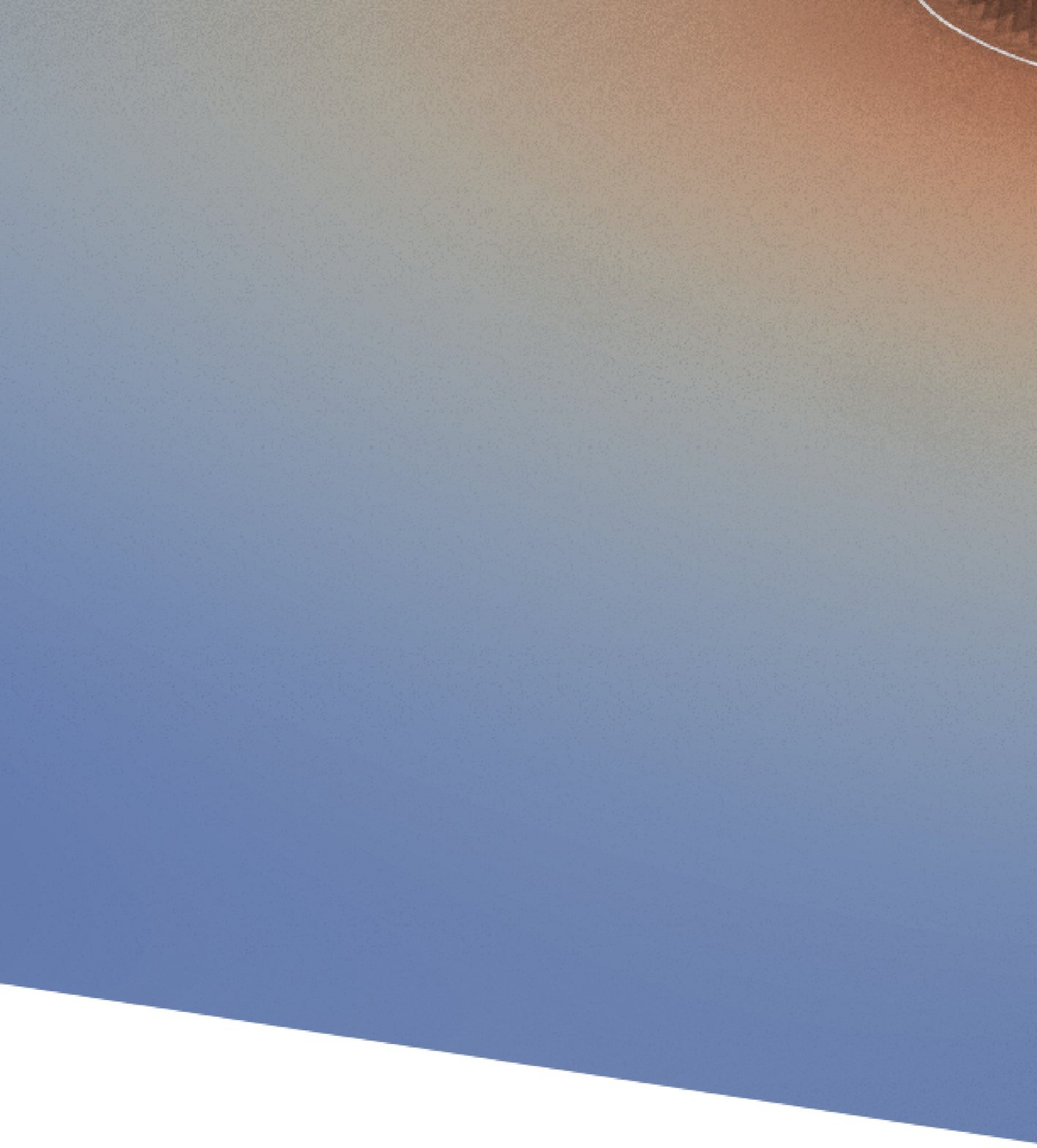}}
\end{tabular}
\caption{\textbf{Illustration of the benchmark problem} considering a surface $\mathcal{S}$ with an isolated ``bump" given by the graph of a function over the two-dimensional domain $\mathrm{\Omega}=[-2,\,2]^2$. \textbf{(a)} Top view of the center part with the initial condition $f_0$ shown by color. The smallest white circle encloses the area where $f_0 \neq 0$ (radius $\sigma = 0.2$). The middle mid-sized white circle encloses the ``bump" (radius $r = 0.25$). The solutions computed by the different numerical methods are probed and compared at the two surface points above $\boldsymbol{x}_0 = (-0.5,\, 0)$ (green dot) and $\boldsymbol{x}_1 = 0.25(-\sqrt{2},\, \sqrt{2})$ (orange dot). The largest white circle of radius $0.5$ and the dotted diagonal line help locate the point $\boldsymbol{x}_1$ with respect to the ``bump".
\textbf{(b)} Perspective view of the surface with the Surface DC-PSE solution at final time $t = 1$ shown by color. For reference, we also show the same circles as shown in (a) for the initial condition and the ``bump", as well as the probe points $\boldsymbol{x}_0$ and $\boldsymbol{x}_1$.}
\label{fig:graph_surf}
\end{figure}

We validate Surface DC-PSE by comparing it with surface Finite Element Methods (FEM) on a test case specifically developed for benchmarking FEM solutions of surface problems \cite{graph_surface_FEM}.
The benchmark problem considers a two-dimensional surface $\mathcal{S}$ with an isolated ``bump'' described by the graph of the function $u(\boldsymbol{x}) = \alpha \zeta(\| \boldsymbol{x} - \boldsymbol{p} \| / r)$, where $\alpha \ge 0$ is the height of the bump, $\boldsymbol{p} = (-0.5,\,0)$ is the position of its center, and $r =0.25$ is the bump radius (see Fig.~\ref{fig:graph_surf}a). The function $\zeta(d) = \exp\big(- \frac{1}{1-d^2}\big)$ for $d < 0.975$ and $0$ otherwise is a cut-off compressed Gaussian. The surface is thus defined as $\mathcal{S} = \{\boldsymbol{x}_{\mathcal{S}} = (x_1, \, x_2, \, u(x_1,x_2)) \in\mathbb{R}^3 \,| \, \boldsymbol{x} = (x_1, \, x_2) \in \mathrm{\Omega}\}$ over the square $\mathrm{\Omega} = [-2,\,2]^2 \subset \mathbb{R}^2$. This surface is asymmetric along the $x$ direction, containing regions of both positive and negative Gauss curvature (see Fig.~\ref{fig:graph_surf}b).

We numerically solve the diffusion equation $\partial_t f (\boldsymbol{x}_{\mathcal{S}},t)= \Lap _{\mathcal{S}} f$ on $\mathcal{S}$ with no-flux Neumann boundary conditions at all borders of the domain $\mathrm{\Omega}$.
The initial condition is the truncated Gaussian $f_0 = f(\boldsymbol{x}_{\mathcal{S}},0) = \sigma^{-2}\zeta(\| \boldsymbol{x}_{\mathcal{S}} \|/\sigma)$ centered at $\boldsymbol{x}_{\mathcal{S}}=(0,0,0)$ with $\sigma = 0.2$ (Fig.~\ref{fig:graph_surf}a).

The collocation points for Surface DC-PSE are regularly distributed in the flat parts of the surface, while in the region of the bump ($[-0.75,-0.25] \times [-0.25,0.25]$) we use the Fibonacci sphere technique \cite{Gonzlez2009} to place the points. This results in a total of $N_p = 18\,439$ points on the surface, of which $16\,441$ are in the flat parts of the surface. The resulting average point spacing is $h = 0.03125$. We use a normal spacing of $\delta n = h$ to generate the surface operators with $N_n=3$ layers to either side of the surface. We choose $r_c = 2.9\delta n$ as the operator support and $r = 2$ as the desired order of convergence. The Neumann boundary conditions at the edge of the surface are imposed using the method of images \cite{cottet_koumoutsakos_2000} with around 3000 ``ghost points'' outside the domain. Time integration over $t \in [0,1]$ is done using the explicit fifth-order Dormand-Prince method~\cite{DORMAND198019} with a time-step size of $\delta t = 10^{-4}$.

We qualitatively compare the solution computed by Surface DC-PSE with those obtained by Finite Element Methods using the data reported in Ref.~\cite{graph_surface_FEM} for Surface FEM (SFEM) \cite{dziuk_finite_2013}, Intrinsic Surface FEM (ISFEM) \cite{bachini_intrinsic_2021}, trace FEM (TraceFEM) \cite{olshanskii_trace_2017}, and diffuse interface (DI) FEM \cite{Lehrenfeld2017,diffuse_interface_2006}. We report the results obtained by high-resolution SFEM and DI FEM. The solutions obtained using a lower-resolution SFEM, as well as using ISFEM and TraceFEM are visually indistinguishable and reported elsewhere~\cite{graph_surface_FEM}.

\begin{figure}[h!]
	\setlength{\tabcolsep}{0pt}
	\centering
	\begin{tabular}{cc}
		\subfloat[]{\includegraphics[width=0.511\textwidth]{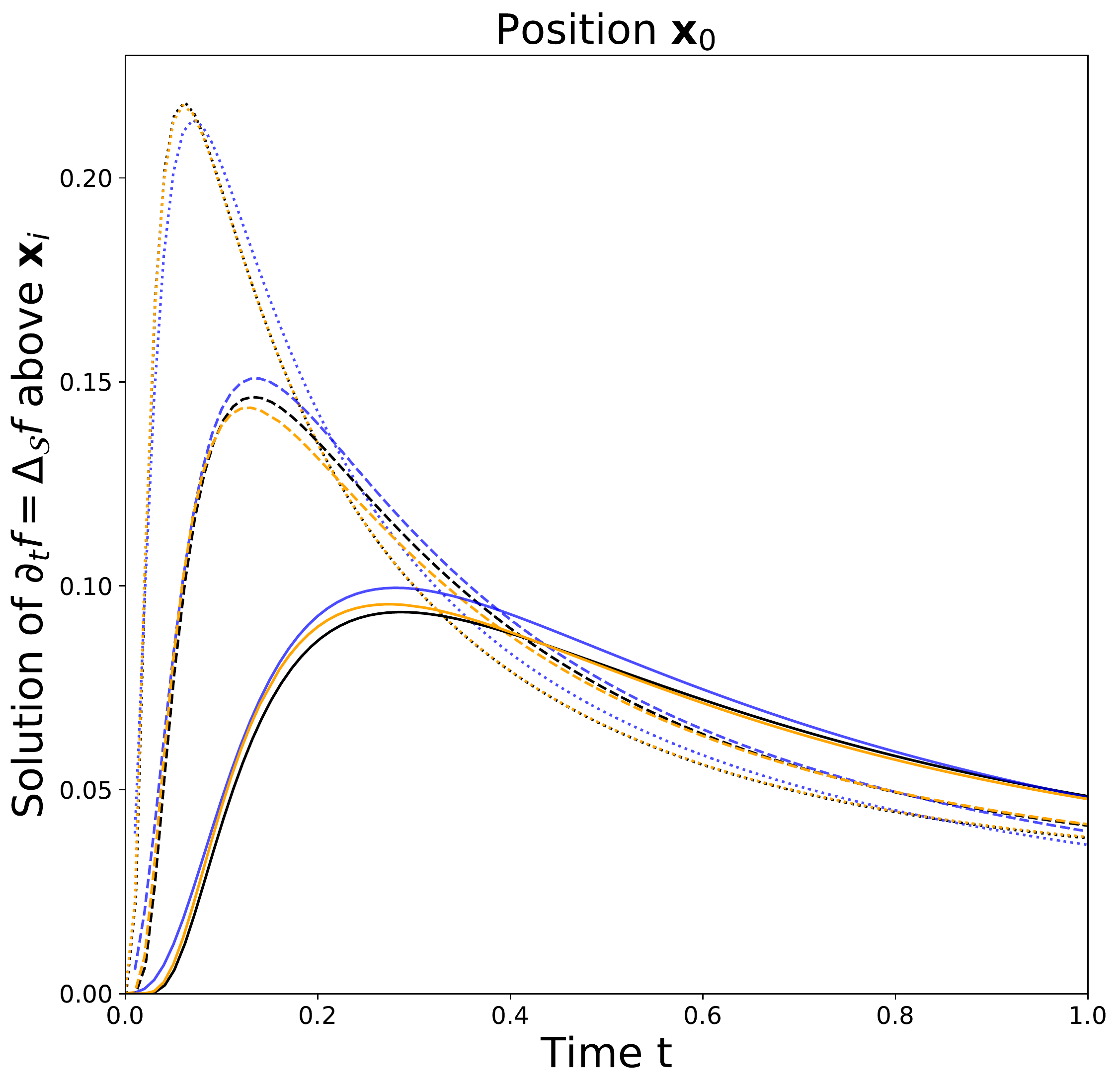}} &
		\subfloat[]{\includegraphics[width=0.49\textwidth]{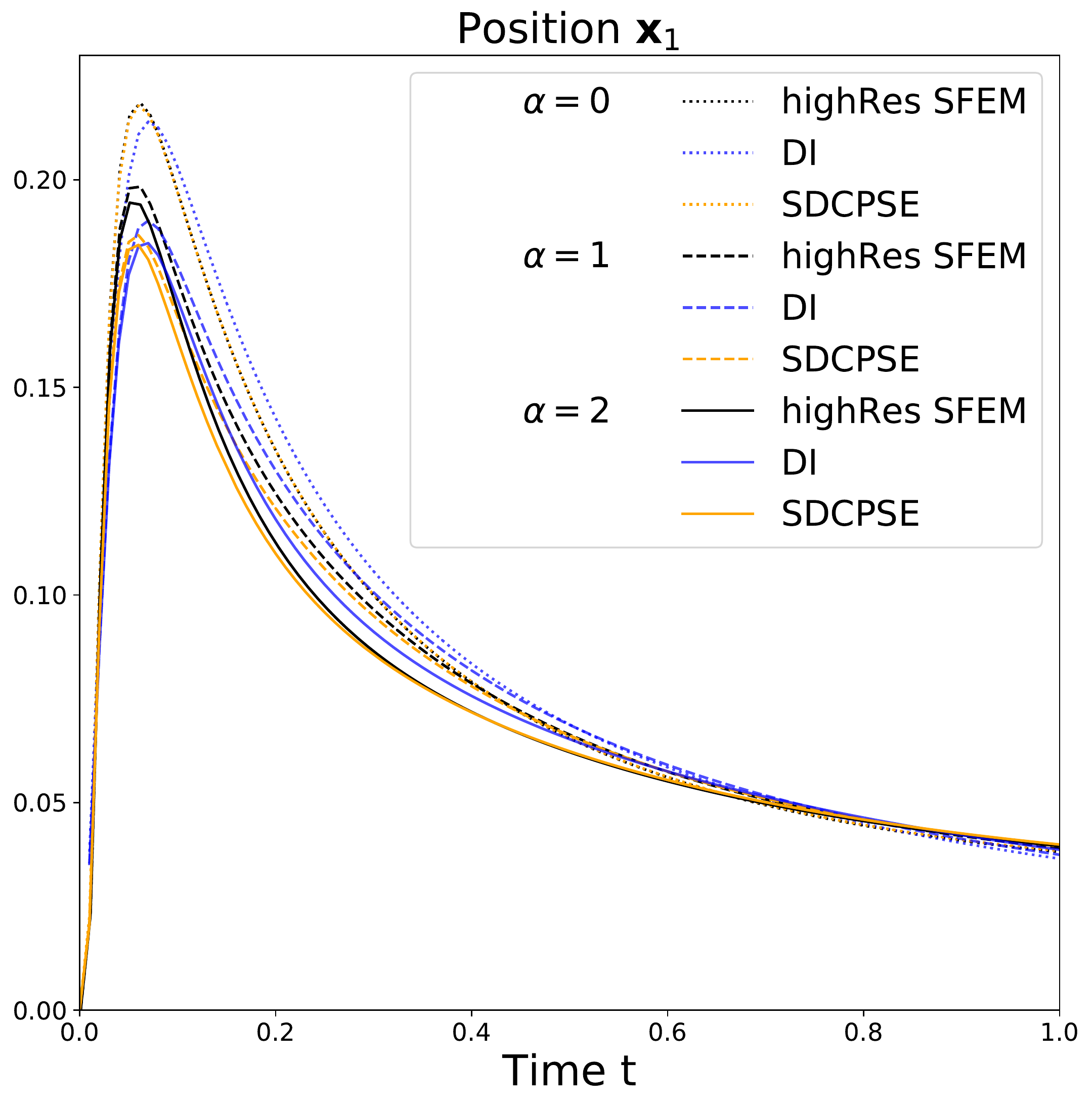}}
	\end{tabular}
	\caption{\textbf{Comparison of Surface DC-PSE (SDCPSE) with two different surface Finite Element Methods (SFEM and DI)} on the  benchmark described in Fig.~\ref{fig:graph_surf} and in the main text, solving the diffusion equation on a graph surface with an isolated  ``bump''~\cite{graph_surface_FEM}. Orange lines correspond to solutions obtained using Surface DC-PSE, blue lines correspond to those from DI FEM, and black lines correspond to those using high-resolution SFEM (reference solution). Solutions for three different bump heights $\alpha$ are displayed by different line styles (see inset legend).
	\textbf{(a)} Plot of the solution $f(\boldsymbol{x}_{\mathcal{S}},t)$ at
	the bump maximum above the point $\boldsymbol{x}_0 = (-0.5,0)$.
	\textbf{(b)} Plot of the solution on the surface above the point $\boldsymbol{x}_1 = 0.25(-\sqrt{2},\sqrt{2})$.}
	\label{fig:graph_surf_results}
\end{figure}

The comparison is done at two points on the surface above $\boldsymbol{x}_0 = \boldsymbol{p} = (-0.5,\,0)$ and $\boldsymbol{x}_1 = 0.25(-\sqrt{2},\, \sqrt{2})$ for three different bump heights $\alpha \in \{0, \, 1, \, 2\}$~\cite{graph_surface_FEM}.
There is no analytical solution for this test case, but a high-resolution SFEM solution serves as a reference. This ``highRes SFEM''
reference solution was obtained using a ten-times finer mesh than Surface DC-PSE, with the finest space resolution (on the bump) set to $h \approx 0.0027$, a time step size of $\delta t = 10^{-4}$, a standard BDF-2 scheme, and a polynomial order of two~\cite{graph_surface_FEM}. Results are shown in Fig.~\ref{fig:graph_surf_results}.

We observe excellent agreement between the Surface DC-PSE solution and the ``highRes SFEM'' reference solution for $\alpha = 0$ (dotted lines in Fig.~\ref{fig:graph_surf_results}). For quantitative comparison, we report the absolute difference in the peak values of the solutions, $e_\text{peak}$. At both probe points, $e_\text{peak} \approx 7 \cdot 10^{-4}$ for the flat case, which is comparable to the other FEM methods \cite{graph_surface_FEM}.

For the curved cases ($\alpha \in \{1, \, 2\}$), Surface DC-PSE is still in very good agreement with the reference solution and closer to it than DI FEM (i.e., the $L_2$ difference between highRes SFEM and Surface DC-PSE is smaller in all cases than between highRes SFEM and DI FEM). For $\alpha = 1$ (dashed lines in Fig.~\ref{fig:graph_surf_results}), the absolute difference in the peak values between Surface DC-PSE and highRes SFEM is $e_\text{peak} \approx 2.6 \cdot 10^{-3}$ above $\boldsymbol{x}_0$ and $e_\text{peak} \approx 1.2 \cdot 10^{-2}$ above $\boldsymbol{x}_1$. The differences for DI FEM are of the same order.
For the highest bump ($\alpha = 2$, solid lines in Fig.~\ref{fig:graph_surf_results}), we observe similar $e_\text{peak}$ values and overall behavior for all methods. Consistently, and as expected, the curves have decreasing peak values for increasing $\alpha$.

\section{Conclusions}\label{sec:Conclusion}
We have presented a meshfree collocation scheme for consistently approximating intrinsic differential operators on smooth curved surfaces represented by point clouds. The present scheme is based on the DC-PSE method for discretizing differential operators on (almost) arbitrarily distributed collocation points in flat spaces~\cite{schrader_discretization_2010}. We have derived the present surface-intrinsic version by realizing that the kernel evaluations can be factored out across points created by {\em exact} constant orthogonal extension, and that the partial sums over the kernels can be precomputed and stored on the surface points only, defining effective surface kernels.
This yields a method that is easy to implement and computationally efficient, as it only requires storing points {\em on} the surface. In this sense, Surface DC-PSE combines features from embedding methods with features from embedding-free methods. The operators are determined in an embedding formulation, but result in a surface-intrinsic algorithm for operator evaluation.

We have verified the method in different test cases with known analytical solutions. This included evaluating the Laplace-Beltrami operator on the unit circle and the unit sphere, solving Poisson equations on the unit circle and the unit sphere using an implicit solver, and computing mean and Gauss curvature of an ellipsoid via approximation of the embedded shape tensor. In all cases, the Surface DC-PSE method converged as expected. We then applied the method to compute the mean curvature of the Stanford bunny, showing an application to a non-parametric surface.
We expect DC-PSE to be more robust (i.e., better conditioned) than Surface MLS in complex geometries \cite{bourantas_using_2016} and likely computationally more efficient, since Surface MLS uses local moving frames in co-dimension 2, which is different from the present tubular approach in co-dimension 1.
Finally, we validated Surface DC-PSE against two different Finite Element Methods (FEM) for surface problems, showing excellent agreement with the reference solution.

Despite its advantages, Surface DC-PSE also has a number of limitations: First, the normal field is required as an input, which can be limiting or introduce additional errors if the normals need to be numerically approximated. Second, for a given point distribution, the numerical error is limited by the curvature of the represented surface and depends on the average spacing between the surface points and the normally extended points (see Figs.~\ref{fig:sdcpse}b and \ref{fig:ellipCurv}b). The required minimum resolution can be determined based on an approximation of the curvature.
Third, Surface DC-PSE is only correct for orientable surfaces that possess a non-intersecting tubular neighborhood (i.e., a tubular network) with a radius at least as large as the kernel radius everywhere. This guarantees that the line segments along which the surface points are extended never intersect. For surfaces that possess any non-intersecting tubular neighborhood, this can always be achieved by choosing the resolution $h$ (locally) sufficiently small, since the tube radius scales with the radius of the DC-PSE kernel.
Fourth, we limited our discussion to surfaces in co-dimension 1. While it may be possible to extend Surface DC-PSE to co-dimension 2 (e.g., curves embedded in $\mathbb{R}^3$), the method is likely not computationally efficient in those cases, as it would require constant orthogonal extension in the whole two-dimensional normal space of each surface point.
Lastly, constructing the Surface DC-PSE kernels is computationally expensive, as it involves solving a small linear system of equations at each surface point in order to determine the embedding-space DC-PSE kernels. For Eulerian simulations, where the collocation points do not move, the kernels have to be determined once at the beginning of the simulation. However, if points move, e.g.~in a Lagrangian simulation or simulations involving deforming surfaces, the kernels need to be recomputed at each time step. While the corresponding cost may be amortized by a gain in accuracy and stability~\cite{schrader_discretization_2010}, it is still significant.

In future work, we will consider extensions of Surface DC-PSE to Lagrangian problems involving moving and deforming surfaces. This also includes simulations of {\em deformable} surfaces, where the surface deformation itself results from intrinsic force-balance equations~\cite{salbreux_mechanics_2017,mietke_self-organized_2019}. We will also consider coupling Surface DC-PSE with regular DC-PSE in the surrounding space in order to describe coupled bulk-surface phenomena, as well as adaptive-resolution Surface DC-PSE with resolution and tube radius locally adjusted to the curvature of the surface in order to maintain error equi-distribution \cite{reboux_self-organizing_2012}.

In summary, we have extended the meshfree collocation method DC-PSE to scalar-valued problems on curved surfaces, requiring only intrinsic surface points. Like DC-PSE, also Surface DC-PSE computes the operator kernels numerically at runtime and is consistent for any desired order of convergence $r$. This makes the presented algorithm particularly attractive for higher-order intrinsic operators, such as the fourth-order operators in Fig.~\ref{fig:SphPois}, and for determining the system matrices of implicit equations on surfaces or implicit time integration schemes.

\section*{acknowledgements}
   We thank Dr.~Simon Praetorius and Prof.~Axel Voigt, both from the Faculty of Mathematics at Technische Universit\"{a}t Dresden, for providing the benchmark data of the FEM methods in Section \ref{sec:fem}. We thank Bryce Palmer (Michigan State University \& Flatiron CCB) for his useful comments on the manuscript.

\section*{Disclosures and Declarations}

\subsection*{Competing interests}
The authors have no competing interests to declare that are relevant to the content of this article.

\subsection*{Funding}
This work was funded by the German Research Foundation (DFG, Deutsche Forschungsgemeinschaft) under grants FOR-3013 (``Vector- and tensor-valued surface PDEs'') (author A.F.), GRK-1907 (``RoSI: role-based software infrastructures'') (author A.S.), and SB-350008342 (``OpenPME'') (author P.I.), and by the Federal Ministry of Education and Research (Bundesministerium f\"{u}r Bildung und Forschung, BMBF) under funding code 031L0160 (project ``SPlaT-DM -- computer simulation platform for topology-driven morphogenesis'') (author I.F.S.).

\section*{Data availability}
The C++ source code of the Surface DC-PSE implementation is available in the numerics module of the open-source scalable scientific computing library OpenFPM: \url{https://github.com/mosaic-group/openfpm_numerics}. All examples shown in this manuscript are contained in the ``example'' folder of the OpenFPM repository.

\end{document}